\documentclass[11pt]{article}%
\usepackage{amsmath}
\usepackage{graphicx}
\usepackage{cite}
\usepackage{subfigure}
\usepackage{multirow}
\usepackage{authblk}
\usepackage{indentfirst}
\usepackage{caption}
\usepackage{algorithm}
\usepackage{algorithmicx}
\usepackage{algpseudocode}
\usepackage{graphicx}
\usepackage{bm}
\usepackage{CJK}
\usepackage{geometry}
\usepackage{float}
\usepackage{colortbl}
\usepackage{textcomp,booktabs}
\usepackage{array}
\usepackage{xcolor}%
\usepackage{amsfonts}%
\usepackage{amssymb}

\newcommand{\eqnb}{\begin{equation}}
\newcommand{\eqne}{\end{equation}}

\newtheorem{The}{Theorem}

\newtheorem{Cor}{Corollary}
\newtheorem{Lem}{Lemma}
\newtheorem{Pro}{Proposition}
\newtheorem{Rem}{Remark}
\definecolor{mygray}{gray}{.9}
\definecolor{mypink}{rgb}{.99,.91,.95}
\definecolor{mycyan}{cmyk}{.3,0,0,0}

\begin{document}

\title{Optimal Energy-Efficient Policies for Data Centers through Sensitivity-Based Optimization}
\author{Jing-Yu Ma\thanks{Jing-Yu Ma and Quan-Lin Li are both with the School of
Economics and Management Sciences, Yanshan University, Qinhuangdao 066004,
China (e-mail: mjy0501@126.com, liquanlin@tsinghua.edu.cn).}, \ Li
Xia\thanks{Li Xia is with the Center for Intelligent and Networked Systems
(CFINS), Department of Automation, TNList, Tsinghua University, Beijing
100084, China (e-mail: xial@tsinghua.edu.cn).}, \ Quan-Lin Li }
\date{}
\maketitle

\begin{abstract}
In this paper, we propose a novel dynamic decision method by
applying the sensitivity-based optimization theory to find the
optimal energy-efficient policy of a data center with two groups of
heterogeneous servers. Servers in Group 1 always work at high energy
consumption, while servers in Group 2 may either work at high energy
consumption or sleep at low energy consumption. An energy-efficient
control policy determines the switch between work and sleep states
of servers in Group 2 in a dynamic way. Since servers in Group 1 are
always working with high priority to jobs, a transfer rule is
proposed to migrate the jobs in Group 2 to idle servers in Group 1.
To find the optimal energy-efficient policy, we set up a
policy-based Poisson equation, and provide explicit expressions for
its unique solution of performance potentials by means of the
RG-factorization. Based on this, we characterize monotonicity and
optimality of the long-run average profit with respect to the
policies under different service prices. We prove that the bang-bang
control is always optimal for this optimization problem, i.e., we
should either keep all servers sleep or turn on the servers such
that the number of working servers equals that of waiting jobs in
Group 2. As an easy adoption of policy forms, we further study the
threshold-type policy and obtain a necessary condition of the
optimal threshold policy. We hope the methodology and results
derived in this paper can shed light to the study of more general
energy-efficient data centers.

\vskip                                                                         0.5cm


\end{abstract}

\section{Introduction}

Data centers have become a core part of the IT infrastructure for Internet
service. Typically, hundreds of thousands of servers are deployed in a data
center to provide ubiquitous computing environments. Tremendous energy
consumption becomes a significant operation expense of data centers. In 2014,
the electricity consumption of data centers in the USA estimated 70 billion
KWh, accounted for 2\% of the national electricity consumption
\cite{Shehabi:2016}. The data centers in the USA are expected to consume
energy 140 TWh and spend \$13 billion energy bills by 2020 \cite{NRDC}, while
these figures in Europe will reach 104 TWh and \$9.6 billion \cite{DeNapoli}.
The energy consumption of data centers consists of three main parts: servers,
networks and cooling, while servers are the major one. It is estimated that
servers consume around 70\% of the total energy consumption in a data center
with tiered architectures \cite{Kliazovich2012}. On the other hand, reducing
the energy consumption of servers also can help reduce the energy consumption
of networking and cooling. Therefore, energy-efficient scheduling of servers
is of significance for the energy management of data centers.

During the last two decades, considerable attention has been paid to
studying the energy efficiency of data centers. An early interesting
observation by Barroso and H\"{o}lzle \cite{Bar:2007} demonstrated
that a lot of data centers were designed to be able to handle the
peak loads effectively, but it directly caused that a significant
number of servers (about 20\%) are often idle in the off-peak
period. Although the idle servers do not provide any service, they
still consume a notable amount of energy. Therefore, it is necessary
to design an energy-efficient mechanism for effectively saving
energy of idle servers. Previous studies demonstrate that a
potential power cutting could be as remarkable as 40\%
\cite{Bod:2012}. For this purpose, a key technique, called an
energy-efficient state `sleep' or `off', was introduced to save
energy for idle servers. See Gandhi et al. \cite{Gan:2012} and Kuehn
and Mashaly \cite{Kue:2015} for more interpretations. In this case,
some queueing models either with server energy-efficient states
(e.g., work, idle, sleep, and off) or with server control policies
(e.g., vacation, setup, and $N$-policy) were developed in the study
of energy-efficient data centers. Queueing theory and Markov (reward
or decision) processes become two useful mathematical tools in
analyzing energy-efficient data centers. See Gandhi \cite{Gan:2013}
and Li et al. \cite{Li:2017} for more details.

Few available studies have applied queueing theory and Markov
processes to performance analysis and optimization of
energy-efficient data centers. Important examples in the recent
literature are remarked as follows. Gandhi et al. \cite{Gan:2010a}
considered a data center with multiple identical servers, the states
of which include work, idle, sleep and off, and their energy
consumption have a decreasing order. One crucial technique given in
Gandhi et al. \cite{Gan:2010a, Gan:2010b} was to develop some
interesting queueing models, for example, the M/M/k queue with setup
times. Since then, some multi-server queues have received attention
(for example, queues with server vacations, queues with either local
setup times or $N$-policy), and they were successfully applied to
energy-efficient management of data centers. Readers may refer to
recent publications for more details, among which are Mazzucco et
al. \cite{Maz:2010a}, Schwartz et al. \cite{Sch:2012}, Gandhi and
Harchol-Balter \cite{Gan:2013a}, Gandhi et al. \cite{Gan:2014},
Maccio and Down \cite{Mac:2015}, Phung-Duc \cite{Phu:2017}, Chen et
al. \cite{Chen:2018}, and Li et al. \cite{Li:2017}.

In the study of energy-efficient data centers, it is a key to
develop effective optimal methods and dynamic control techniques in
data centers. So far, there have been two classes of optimal methods
applied to the analysis of energy-efficient data centers. The first
class is regarded as `static optimization' with two basic steps.
Step one is to set up performance cost (i.e., a suitable
performance-energy tradeoff) of a data center, where the performance
cost can be expressed by means of queueing indexes of the data
center. Step two is to optimize the performance cost with respect to
some key parameters of the data center by using, such as, linear
programming, nonlinear programming, integer programming, and bilevel
programming. The second class is viewed as `dynamic optimization' in
which Markov decision processes or stochastic network optimization
are applied to energy-efficient management of data centers, e.g.,
see Benini et al. \cite{Ben:2000} and Yao et al. \cite{Yao:2014} for
more details.

For the static optimization, some available works have been successfully
conducted according to two key points: The first key point is to emphasize how
to construct a suitable utility function for the performance-energy tradeoff,
which needs to synchronously optimize several different performance measures,
for example, reducing energy consumption, reducing system response time, and
improving quality of service. The second key point is to minimize performance
cost with respect to some crucial parameters of data centers by means of, such
as linear programming and nonlinear programming. On such a research line,
Gandhi et al. \cite{Gan:2010a} recalled two classes of performance-energy
tradeoffs: \textbf{(a) ERWS}, the weighted sum $\beta_{1}[R]+\beta_{2}[E]$ of
the mean response time $[R]$ and the mean power cost $[E]$, where $\beta
_{1},\beta_{2}\geq0$ are weighted coefficients; and \textbf{(b) ERP}, the
product $[R][E]$ of the mean response time and the mean power cost. For the
ERP, Gandhi et al. \cite{Gan:2010a} first described the data center as a queue
to compute the two mean values $[R]$ and $[E]$, and then provided an
optimization method to minimize the ERP. Also, they further analyzed
optimality or near-optimality of several different energy-efficient policies.
In addition, Gandhi \cite{Gan:2013} gave some extended results and a
systematical summarization with respect to minimizing the ERP. Maccio and Down
\cite{Mac:2015} generalized the ERP by Gandhi \cite{Gan:2010a} to a more
general performance cost function as follows.
\[
f\left(  \beta,w\right)  =\sum_{i=1}^{M}\beta_{i}\left(  [R]\right)
^{w_{R,i}}\left(  [E]\right)  ^{w_{E,i}}\left(  [C]\right)  ^{w_{C,i}},
\]
where $[C]$ is the expected cycle rate, and $\beta_{i}$, $w_{R,i}$, $w_{E,i}$
and $w_{C,i}$ for $1\leq i\leq M$ are nonnegative weighted coefficients,
$\beta=\left(  \beta_{i}:1\leq i\leq M\right)  $ and $w=\left(  w_{R,i}%
,w_{E,i},w_{C,i}:1\leq i\leq M\right)  $. They used the queueing models to
compute the three mean values $[R]$, $[E]$ and $[C]$, and then provided some
discussion on the optimality of cost function $f\left(  \beta,w\right)  $.
Gebrehiwot et al. \cite{Geb:2016} made another interesting generalization of
the ERP and ERWS by Gandhi \cite{Gan:2010a} through introducing the multiple
intermediate sleep states. Under more general assumptions with general service
and setup times, they computed the two mean values $[R]$ and $[E]$ by means of
some queueing insensitivity, and then discussed the optimality of the ERP and
ERWS. Further, Gebrehiwot et al. \cite{Geb:2016a, Geb:2017} generalized the
FCFS queueing results of the data center with multiple intermediate sleep
states to the processor-sharing discipline and the shortest remaining
processing time (SRPT) discipline, respectively. Different from the ERP and
ERWS, Mitrani \cite{Mit:2011, Mit:2013} considered a data center of $N$
identical servers whose first part contains $m$ servers. The idle or work
state of servers is controlled by two different thresholds: an up threshold
$U$ and a down threshold $D$. He designed a simple three-layer queue to
describe the energy-efficient data center in terms of a new performance cost:
$C=c_{1}L+c_{2}S$, where $L$ and $S$ are the average numbers of jobs present
and of energy consumption, respectively. He provided expressions for computing
the average numbers $L$ and $S$ such that the performance cost $C$ can be
optimized with respect to the three parameters $m$, $U$ and $D$.

However, for the dynamic optimization, little work has been done on
applying Markov decision processes to set up optimal dynamic control
policies for energy-efficient data centers. In general, such a study
is more interesting, difficult and challenging due to the fact that
a complicated queueing model with synchronously multiple control
objectives (e.g., reducing energy consumption, reducing system
response time and guaranteeing quality of service) needs to be
synthetically established in a Markov decision process. For a data
center with multiple identical servers, Kamitsos et al.
\cite{Kam:2010, Kam:2012, Kam:2017} constructed a discrete-time
Markov decision process by uniformization and proved that the
optimal sleep energy-efficient policy is simply hysteretic. Hence,
this problem has a double threshold structure by means of the
optimal hysteretic policy given in Hipp and Holzbaur \cite{Hip:1988}
and Lu and Serfozo \cite{Lu:1984}. On the other hand, as some close
research to energy-efficient data centers, it is worthwhile to note
that the policy optimization and dynamic power management for
electronic systems or equipments were developed well by means of
Markov decision processes and stochastic network optimization.
Important examples include: (a) Discrete-time Markov decision
processes by Benini et al. \cite{Ben:2000} and Yang et al.
\cite{Yan:2017}; (b) Continuous-time Markov decision processes by
Qiu and Pedram \cite{Qiu:1999} and Qiu et al. \cite{Qiu:2001}; (c)
Stochastic network optimization by Yao et al. \cite{Yao:2014} and
Huang and Neely \cite{Hua:2013}; (d) It has become increasingly
important to simplify the method of Markov decision processes such
that more complicated stochastic networks can be analyzed
effectively. On this research line, event-driven techniques of
Markov decision processes have received high attention for the past
one decade. Important examples include the event-driven power
management by \v{S}imuni\'{c} et al. \cite{Sim:2001}, and the
event-driven optimization techniques by Becker et al.
\cite{Bec:2004}, Cao \cite{Cao:2007}, Koole \cite{Koo:1998}, Engel
and Etzion \cite{Eng:2011}, and Xia et al. \cite{Xia:2014}.

The purpose of this paper is to apply the Markov decision processes to set up
an optimal dynamic control policy for energy-efficient data centers. To do
this, we first apply the sensitivity-based optimization theory in the study of
data centers. Note that the sensitivity-based optimization is greatly refined
from the Markov decision processes through re-expressing the Poisson equation
(corresponding to the Bellman optimality equation) by means of several novel
tools, for instance, performance potential and performance difference (see
Cao's book \cite{Cao:2007}). Also, the sensitivity-based optimization theory
can be effectively related to the Markov reward processes (e.g., see Li
\cite{Li:2010} and Li and Cao \cite{Li:2004}) so that it is an effective
dynamic decision method for performance optimization of many practical Markov
systems. The key idea in the sensitivity-based optimization theory is a
performance difference equation that can quantify the performance difference
of a Markov system under any two different policies. The difference equation
gives a clear relation that explains how the system performance is varying
with respect to policies. See an excellent book by Cao \cite{Cao:2007} for
more details. So far, the sensitivity-based optimization theory has been
applied to performance optimization of queueing systems (or networks).
Important examples include an early invited overview by Xia and Cao
\cite{Xia:2012}; the MAP/M/1 queue by Xia et al. \cite{Xia:2017a}; the closed
queueing networks by Xia and Shihada \cite{Xia:2013}, Xia \cite{Xia:2014a} and
Xia and Jia \cite{Xia:2015}; and the open queueing networks by Xia
\cite{Xia:2014b} and Xia and Chen \cite{Xia:2018}. In addition, the
sensitivity-based optimization theory was also applied to network energy
management, for example, the multi-hop wireless networks by Xia and Shihada
\cite{Xia:2015a} and the tandem queues with power constraints by Xia et al.
\cite{Xia:2017b}.

The main contributions of this paper are twofold. The first one is to apply
the sensitivity-based optimization theory to study the optimal
energy-efficient policies of data centers for the first time, in which we
propose a job transfer rule among the server groups such that the sleep
energy-efficient mechanism becomes more effective. Different from previous
works in the literature for applying an ordinary Markov decision process to
dynamic control of data centers, we propose and develop an easier and more
convenient dynamic decision method: sensitivity-based optimization, in the
study of energy-efficient data centers. Crucially, this sensitivity-based
optimization method may open a new avenue to the optimal energy-efficient
policy of more complicated data centers. The second contribution is to
characterize the optimal energy-efficient policy of data centers. We set up a
policy-based Poisson equation and provide explicit expression for its unique
solution by means of the RG-factorization. Based on this, we analyze the
monotonicity and optimality of the long-run average profit with respect to the
energy-efficient policies under some restrained service prices. We obtain the
structure of optimal energy-efficient policy. Specifically, we prove that the
bang-bang control is optimal for this problem, which significantly reduces the
large search space. We also provide an effective way to design and verify the
threshold-type mechanism in practice, which is of great significance to solve
the mechanism design problem of energy-efficient data centers. Therefore, the
results of this paper give new insights on understanding not only mechanism
design of energy-efficient data centers, but also applying the
sensitivity-based optimization to dynamic control of data centers. We hope
that the methodology and results given in this paper can shed light to the
study of more general energy-efficient data centers.

The remainder of this paper is organized as follows. In Section 2, we describe
the problem of an energy-efficient data center with two groups of different
servers. In Section 3, for the energy-efficient data center, we first
establish a policy-based continuous-time birth-death process with finite
states. Then we define a suitable reward function with respect to states and
policies of the birth-death process. Based on this, we formulate a dynamic
optimization problem to find the optimal energy-efficient policy of the data
center. In Section 4, we set up a policy-based Poisson equation and provide
explicit expression for its unique solution by means of the RG-factorization.
In Section 5, we define a perturbation realization factor of the policy-based
control process of the data center, and analyze how the service price impacts
on the perturbation realization factor. In Section 6, we use the Poisson
equation to derive a useful performance difference equation. Based on this, we
discuss the monotonicity and optimality of the long-run average profit with
respect to the energy-efficient policies, and prove the optimality of the
bang-bang control. In Section 7, we use the Poisson equation to further study
a class of threshold energy-efficient policies, and obtain the necessary
condition of the optimal threshold policy. Finally, we give some discussions
and conclude this paper in Section 8.

\section{Problem Description}

In this section, we give a problem description of the energy-efficient problem
in data centers. As the large variation of working loads in data centers, it
is widely adopted to organize and operate the data center in a multiple-tier
architecture such that the on/off scheduling can be performed in different
tiers to save energy \cite{Kliazovich2012}. In this paper, we study a data
center with two groups of heterogeneous servers. There is no waiting room for
jobs in the data center, which can be viewed as a loss queue. The assumption
of loss queue model is reasonable for data centers and it is also widely used
in telephone systems, computer networks, cloud computing, and so on
\cite{Has:2015,Hon:2013,Tan:2012}. In what follows we provide a detailed
problem description for the data center.

\textbf{Server groups:} The data center contains two server groups: Groups 1
and 2, each of which is also one of the interactive subsystems of the data
center. Groups 1 and 2 contain $n$ and $m$ servers, respectively. Hence the
data center contains $n+m$ servers. Servers in the same group are homogeneous
and in different groups are heterogeneous. Note that Group 1 is viewed as a
base-line group whose servers are always at the work state to guarantee a
necessary service capacity in the data center. Each server in Group 1 consumes
an amount of energy per unit of time. By contrast, Group 2 is regarded as a
reserved group whose servers may either work or sleep so that each of the $m$
servers can switch its state between work and sleep. If one server in Group 2
is at the sleep state, then it consumes a smaller amount of energy to maintain
the sleep state.

\textbf{Power consumption: }The power consumption rates (i.e., power
consumption per unit of\ time) for the two groups of servers are described as:
$P_{1,W}$ and $P_{2,W}$ for the work state in Groups 1 and 2, respectively;
and $P_{2,S}$ only for the sleep state in Group 2. We assume that
$0<P_{2,S}<P_{2,W}$.

\textbf{Arrival processes:} The arrivals of jobs at the data center are a
Poisson process with arrival rate $\lambda$. Each arriving job is assigned to
a server of the two groups according to the following allocation rules:

(a) \textit{Each server in Group 1 must be fully utilized so that Group 1
provides a priority service over Group 2}. If Group 1 has some idle servers,
then the arriving job immediately enters the idle server in Group 1 and
receives service there. Furthermore, if all the servers of Group 1 are busy
but Group 2 has some idle servers, then the arriving job immediately enters an
idle server in Group 2 and receives service there.

(b) \textit{No waiting room}. A job can be served at an idle server or wait at
a sleeping server in Group 2. If all the servers of Groups 1 and 2 are
occupied, then any arriving job must be lost immediately. Note that each
server may contain only one job, hence the total number of jobs in the data
center cannot exceed the number $n+m$.

(c) \textit{Opportunity cost}. Once the data center contains $n+m$ jobs, then
any arriving job has to be lost immediately due to no waiting room. This leads
to an opportunity cost with respect to the job loss.

\textbf{Service processes:} In Groups 1 and 2, the service times provided by
each server are independent and exponential with the service rate $\mu_{1}$
and $\mu_{2}$, respectively. We assume that $\mu_{1}\geq\mu_{2}$ as a
\textit{fast condition}, which makes the prior use of servers in Group~1.

Once a job enters the data center to receive or wait for service, it has to
pay a holding cost in Group 1 or Group 2. We assume a so-called \textit{cheap
condition} that the holding cost in Group 1 is always cheaper than that in
Group 2. The fast and cheap conditions are intuitive to guarantee the prior
use of servers in Group 1. That is, the servers in Group 1 are not only faster
but also cheaper than those in Group 2.

If a job finishes its service at a server and leaves the system, then the data
center can obtain a fixed service revenue from each served job. The service
discipline of each server in the data center is First Come First Serve (FCFS).

\textbf{Transfer rule: }Based on prior use of servers in Group 1, whenever a
server in Group 1 becomes idle, an incomplete-service job (if exists) in Group
2 should be transferred to the idle server in Group 1 to save processing time.
When a job is transferred from Group 2 to Group 1, the data center needs to
pay a transfer cost.

\textbf{Independence: }We assume that all the random variables defined above
in the data center are independent.

\begin{figure}[th]
\centering            \includegraphics[width=12cm]{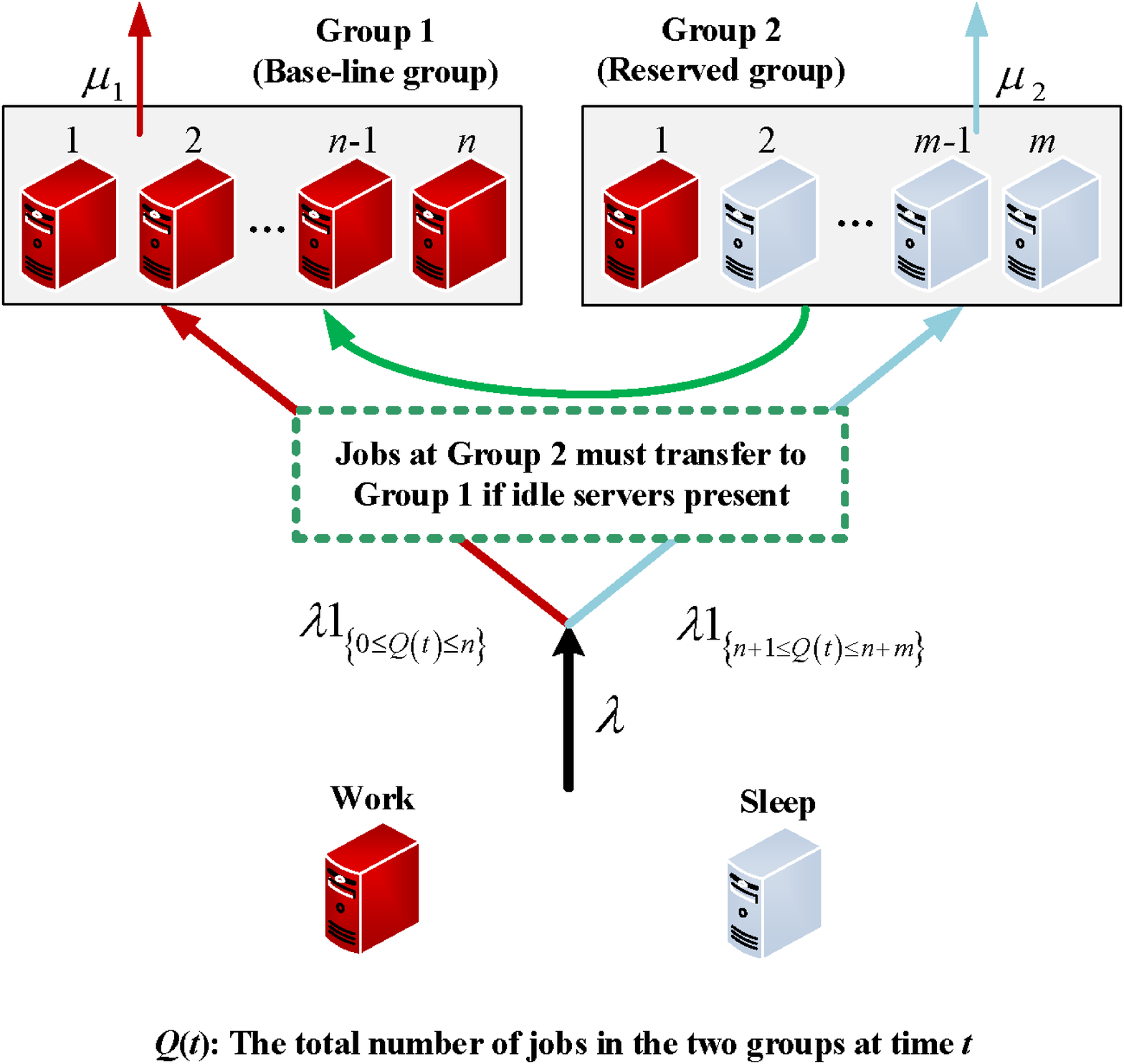}
\caption{Energy-efficient management of data centers.}%
\label{figure:fig-1}%
\end{figure}

Finally, the data center, together with its operational mode and mathematical
notations, are simply depicted in Figure~\ref{figure:fig-1}.

\begin{Rem}
A further interpretation for the transfer rule: When a job is being served at
a server of Group 2, it can be transferred to one idle server of Group 1 and
restart its service when an idle server in Group 1 is available. Note that
each server in Group 1 is not only faster but also cheaper than that in Group
2, that is, the fast and cheap conditions guarantee that servers in Group 1
have priority than those in Group 2. From the memoryless property of
exponential distributions, the restarting service in Group 1 is still faster
and cheaper than its original service in Group 2. Therefore, this transfer
rule effectively supports energy-efficient management of the data center due
to the fact that the servers of Group 1 are fully utilized while servers of
Group 2 are kept in sleep state as many as possible.
\end{Rem}

\begin{Rem}
Although some authors (e.g., see Gandhi et al. \cite{Gan:2010a}, Mitrani
\cite{Mit:2011, Mit:2013} and Maccio and Down \cite{Mac:2015}) analyzed the
energy-efficient data center with two groups of servers, where one is the
base-line group and the other is the reserved or subsidiary group, all the
servers in the two groups given in their papers are assumed to be identical.
From such a point of view, it is easy to understand that those queueing models
given in their works are the same as only Group 2 of our paper. On the other
hand, it is noted that for our queueing model, Group 1 has a large influence
on analysis of Group 2 due to introducing new \textbf{transfer rules}. In
fact, our queueing model here has been discussed in Li et al. \cite{Li:2017}
with a more general setting.
\end{Rem}

\section{Optimization Model Formulation}

In this section, for the energy-efficient data center, we first establish a
policy-based continuous-time Markov process, and show that its infinitesimal
generator has the simple structure of a birth-death process with finite
states. Then, we define a suitable reward function with respect to both states
and policies of the birth-death process. Based on this, we set up a dynamic
optimization model to deal with the optimal energy-efficient policy of the
data center.

In the data center with Group 1 of $n$ servers and Group 2 of $m$ servers, we
need to introduce both `states' and `policies' to express stochastic dynamics
of this data center. Let $I\left(  t\right)  $ and $J\left(  t\right)  $ be
the number of jobs in Groups 1 and 2, respectively. Then, $\left(  I\left(
t\right)  ,J\left(  t\right)  \right)  $ is regarded as the state of the data
center at time $t$. Let all the cases of state $\left(  I\left(  t\right)
,J\left(  t\right)  \right)  $ form a state space as follows.%
\[
\bm      \Omega=\left\{  \left(  0,0\right)  ,\left(  1,0\right)
,\ldots,\left(  n,0\right)  ,\left(  n,1\right)  ,\ldots,\left(  n,m\right)
\right\}  .
\]
Note that such a state $\left(  n-i,j\right)  $ for $1\leq i\leq n$ and $1\leq
j\leq m$ does not exist according to the transfer rule. However, the policies
are defined with a little bit more complication. Let $d_{i,j}$ be the number
of servers turned on in Group 2 at the state $\left(  I\left(  t\right)
,J\left(  t\right)  \right)  =\left(  i,j\right)  $, for $i=0,1,\ldots,n$ and
$j=0,1,\ldots,m$. From the problem description in Section 2, it is easy to see
that%
\begin{equation}
d_{i,j}=\left\{
\begin{array}
[c]{ll}%
0, & i=0,1,\ldots,n-1,j=0,\\
0, & i=n,j=0,\\
0,1,\ldots,m, & i=n,j=1,2,\ldots,m.
\end{array}
\right.  \label{d-1}%
\end{equation}
Now, we provide an interpretation for the above expression: If $I\left(
t\right)  =0,1,\ldots,n-1$, then $J\left(  t\right)  =0$ due to the transfer
rule. In this case with $i=0,1,\ldots,n$, $d_{i,0}$ is taken as zero due to
the energy-efficient cause. While once $I\left(  t\right)  =n$, $J\left(
t\right)  $ may be any element in the set $\left\{  0,1,\ldots,m\right\}  $.
In this case with $i=n$ and $j\neq0$, $d_{n,j}$ may be taken as any element in
the set $\left\{  0,1,\ldots,m\right\}  $.

Corresponding to each case of state $\left(  I\left(  t\right)  ,J\left(
t\right)  \right)  $ at time $t$, we define a time-homogeneous policy as
\[
\bm  d=\left(  d_{0,0},d_{1,0},\ldots,d_{n,0};d_{n,1},d_{n,2},\ldots
,d_{n,m}\right)  .
\]
It follows from (\ref{d-1}) that
\begin{equation}
\bm  d=\left(  0,0,\ldots,0;d_{n,1},d_{n,2},\ldots,d_{n,m}\right)  .
\label{d-2}%
\end{equation}
Let all the possible policies $\bm  d$ given in (\ref{d-2})\ compose a policy
space as follows.
\[
\mathcal{D}=\left\{  \bm  d:\bm  d=\left(  0,0,\ldots,0;d_{n,1},d_{n,2}%
,\ldots,d_{n,m}\right)  ,d_{i,j}\in\left\{  0,1,\ldots,m\right\}  ,\text{
}\left(  i,j\right)  \in\bm  \Omega\right\}  .
\]
Let $\mathbf{X}^{\left(  \bm  d\right)  }(t)=\left(  I\left(  t\right)
,J\left(  t\right)  \right)  ^{(\bm  d)}$ be the system state at time $t$
under any given policy $\bm  d\in\mathcal{D}$. Then, $\{\mathbf{X}^{\left(
\bm  d\right)  }(t):t\geq0\}$ is a policy-based continuous-time birth-death
process on the state space $\bm  \Omega$ whose infinitesimal generator is
given by%
\begin{equation}
\mathbf{B}^{\left(  \bm  d\right)  }=\left(
\begin{array}
[c]{cccccc}%
-\lambda & \lambda &  &  &  & \\
\mu_{1} & -\left(  \lambda+\mu_{1}\right)  & \lambda &  &  & \\
& \ddots\text{ \ \ \ \ \ \ } & \ddots\text{ \ \ \ \ \ \ \ } & \ddots\text{
\ \ \ \ \ \ } &  & \\
& n\mu_{1} & -\left(  \lambda+n\mu_{1}\right)  & \lambda &  & \\
&  & \nu\left(  d_{n,1}\right)  & -\left[  \lambda+\nu\left(  d_{n,1}\right)
\right]  & \lambda & \\
&  & \text{\ \ \ \ \ \ \ }\ddots & \ \ \ \ \ \ \  & \ddots\text{
\ \ \ \ \ \ \ \ \ \ }\ddots & \text{ \ \ \ \ \ }\\
&  &  & \nu\left(  d_{n,m-1}\right)  & -\left[  \lambda+\nu\left(
d_{n,m-1}\right)  \right]  & \lambda\\
&  &  &  & \nu\left(  d_{n,m}\right)  & -\nu\left(  d_{n,m}\right)
\end{array}
\right)  , \label{eq-60}%
\end{equation}
where $\nu\left(  d_{n,j}\right)  =n\mu_{1}+\left(  d_{n,j}\wedge j\right)
\mu_{2}$ for $j=1,2,\ldots,m$, and $a\wedge b$ denotes the minimal one between
two real numbers $a$ and $b$. Note that $d_{n,j}\geq0$, it is clear that
$\nu\left(  d_{n,j}\right)  \geq n\mu_{1}>0$. Thus, the birth-death process
$\mathbf{B}^{\left(  \bm  d\right)  }$ must be irreducible, aperiodic and
positive recurrent for any given policy $\bm  d\in\mathcal{D}$. In this case,
we write the stationary probability vector of the Markov process $\left\{
\mathbf{X}^{\left(  \bm  d\right)  }(t):t\geq0\right\}  $ under a policy $\bm
d\in\mathcal{D}$.
\begin{equation}
\bm  \pi^{\left(  \bm  d\right)  }=\left(  \pi^{\left(  \bm  d\right)
}\left(  0,0\right)  ,\pi^{\left(  \bm  d\right)  }\left(  1,0\right)
,\ldots,\pi^{\left(  \bm  d\right)  }\left(  n,0\right)  ,\pi^{\left(  \bm
d\right)  }\left(  n,1\right)  ,\ldots,\pi^{\left(  \bm  d\right)  }\left(
n,m\right)  \right)  . \label{eq-1-1}%
\end{equation}

Obviously, the stationary probability vector $\bm      \pi^{\left(  \bm
d\right)  }$ is the unique solution to the system of linear equations: $\bm
\pi^{\left(  \bm      d\right)  }\mathbf{B}^{\left(  \bm      d\right)  }=\bm
0$ and $\bm      \pi^{\left(  \bm      d\right)  }\bm      e=1$, where $\bm
e$ is a column vector of ones with proper dimension. We write
\begin{equation}
\xi_{i,0}=\dfrac{\lambda^{i}}{i!\mu_{1}^{i}},\text{ \ }i=0,1,\ldots,n,
\label{eq-1-01}%
\end{equation}
and%
\begin{equation}
\xi_{n,j}^{\left(  \bm  d\right)  }=\dfrac{\lambda^{n}}{n!\mu_{1}^{n}}%
\dfrac{\lambda^{j}}{\underset{i=1}{\overset{j}{\Pi}}\nu\left(  d_{n,i}\right)
},\text{ \ }j=1,2,\ldots,m, \label{eq-1-02}%
\end{equation}%
\begin{equation}
b^{\left(  \bm  d\right)  }=\overset{n}{\underset{i=0}{\sum}}\xi
_{i,0}+\overset{m}{\underset{j=1}{\sum}}\xi_{n,j}^{\left(  \bm  d\right)  }.
\label{eq-1-4}%
\end{equation}
It follows from Subsection 1.1.4 of Chapter 1 in Li \cite{Li:2010} that
\begin{align}
\pi^{\left(  \bm  d\right)  }\left(  i,0\right)   &  =\frac{1}{b^{\left(  \bm
d\right)  }}\xi_{i,0},\text{ \ }i=0,1,\ldots,n;\label{eq-1-2}\\
\pi^{\left(  \bm  d\right)  }\left(  n,j\right)   &  =\frac{1}{b^{\left(  \bm
d\right)  }}\xi_{n,j}^{\left(  \bm  d\right)  },\text{ \ }j=1,2,\ldots,m.
\label{eq-1-3}%
\end{align}

For any two vectors $\bm{a}  =\left(  a_{1},a_{2},\ldots,a_{K}\right)  $ and
$\bm{b}  =\left(  b_{1},b_{2},\ldots,b_{K}\right)  $, we say that $\bm{a}
\geq$ $\bm{b}  $ if $a_{i}\geq b_{i}$ for any $1\leq i\leq K$. The following
proposition provides an interesting observation on how a policy $\bm
d\in\mathcal{D}$ influences the stationary probability vector $\bm
\pi^{\left(  \bm      d\right)  }$.

\begin{Pro}
For any two given policies $\bm                                 d_{1},\bm
d_{2}\in\mathcal{D}$ with $\bm                               d_{1}\geq\bm
d_{2}$, then
\[
\pi^{\left(  \bm                     d_{1}\right)  }\left(  i,0\right)
\geq\pi^{\left(  \bm                d_{2}\right)  }\left(  i,0\right)  ,\text{
\ }i=0,1,\ldots,n.
\]
\end{Pro}

\textbf{Proof:} For any two given policies $\bm    d_{1},\bm    d_{2}%
\in\mathcal{D}$ with $\bm    d_{1}\geq\bm    d_{2}$, then it follows from
(\ref{eq-1-01}) that%
\[
\xi_{i,0}^{\left(  \bm    d_{1}\right)  }=\xi_{i,0}^{\left(  \bm
d_{2}\right)  }=\dfrac{\lambda^{i}}{i!\mu_{1}^{i}},\text{ \ }i=0,1,\ldots,n.
\]
If $\bm    d_{1}\geq\bm    d_{2}$, then for each $j=1,2,\ldots,m$, it is clear
that $d_{n,j}^{1}\geq d_{n,j}^{2}$, this gives%
\[
\nu\left(  d_{n,j}^{1}\right)  =n\mu_{1}+\left(  d_{n,j}^{1}\wedge j\right)
\mu_{2}\geq n\mu_{1}+\left(  d_{n,j}^{2}\wedge j\right)  \mu_{2}=\nu\left(
d_{n,j}^{2}\right)  ,
\]
hence it follows from (\ref{eq-1-02}) that%
\[
\xi_{n,j}^{\left(  \bm    d_{1}\right)  }=\dfrac{\lambda^{n}}{n!\mu_{1}^{n}%
}\dfrac{\lambda^{j}}{\underset{i=1}{\overset{j}{\Pi}}\nu\left(  d_{n,i}%
^{1}\right)  }\leq\dfrac{\lambda^{n}}{n!\mu_{1}^{n}}\dfrac{\lambda^{j}%
}{\underset{i=1}{\overset{j}{\Pi}}\nu\left(  d_{n,i}^{2}\right)  }=\xi
_{n,j}^{\left(  \bm    d_{2}\right)  }.
\]
It is easy to see from (\ref{eq-1-4}) that%
\[
b^{\left(  \bm    d_{1}\right)  }=\overset{n}{\underset{i=0}{\sum}}\xi
_{i,0}+\overset{m}{\underset{j=1}{\sum}}\xi_{n,j}^{\left(  \bm
d_{1}\right)  }\leq\overset{n}{\underset{i=0}{\sum}}\xi_{i,0}+\overset
{m}{\underset{j=1}{\sum}}\xi_{n,j}^{\left(  \bm    d_{2}\right)  }=b^{\left(
\bm  d_{2}\right)  }.
\]
Thus, it follows (\ref{eq-1-2}) that for each $i=0,1,\ldots,n$,%
\[
\pi^{\left(  \bm    d_{1}\right)  }\left(  i,0\right)  =\frac{1}{b^{\left(
\bm   d_{1}\right)  }}\xi_{i,0}\geq\frac{1}{b^{\left(  \bm    d_{2}\right)  }%
}\xi_{i,0}=\pi^{\left(  \bm    d_{2}\right)  }\left(  i,0\right)  .
\]
This completes the proof. \textbf{{\rule{0.08in}{0.08in}}}

The following theorem provides a useful observation that some special policies
$\bm      d\in\mathcal{D}$ have no effect on both the infinitesimal generator
$\mathbf{B}^{\left(  \bm      d\right)  }$ and the stationary probability
vector $\bm      \pi^{\left(  \bm      d\right)  }$. Note that this theorem
will be necessary and useful for the analysis of policy monotonicity and
optimality in our later study, for example, the proof of Theorem 3.

\begin{The}
\label{The:pi}Suppose that two policies $\bm                     d_{1},\bm
d_{2}\in\mathcal{D}$ satisfy the following two conditions: For each
$j=1,2,\ldots,m$, (a) if $d_{n,j}^{1}\in\left\{  1,2,\ldots,j-1\right\}  $,
then we take $d_{n,j}^{2}=d_{n,j}^{1}$; and (b) if $d_{n,j}^{1}\in\left\{
j,j+1,\ldots,m\right\}  $, we take $d_{n,j}^{2}$ as any element of the set
$\left\{  j,j+1,\ldots,m\right\}  $. We have%
\[
\mathbf{B}^{\left(  \bm                           d_{1}\right)  }%
=\mathbf{B}^{\left(  \bm             d_{2}\right)  },\text{ \ }\bm
\pi^{\left(  \bm   d_{1}\right)  }=\bm                         \pi^{\left(
\bm       d_{2}\right)  }.
\]
\end{The}

\textbf{Proof:} It is easy to see from (\ref{eq-60}) that the first $n+1$ rows
of the matrix $\mathbf{B}^{\left(  \bm                         d_{1}\right)
}$ is the same as those of the matrix $\mathbf{B}^{\left(  \bm
d_{2}\right)  }$.

In what follows we compare the latter $m$ rows of the matrix $\mathbf{B}%
^{\left(  \bm    d_{1}\right)  }$ with those of the matrix $\mathbf{B}%
^{\left(  \bm    d_{2}\right)  }$. For the two policies $\bm    d_{1},\bm
d_{2}\in\mathcal{D}$ satisfying the two conditions (a) and (b), by using
$\nu\left(  d_{n,j}\right)  =n\mu_{1}+\left(  d_{n,j}\wedge j\right)  \mu_{2}%
$, it is clear that for $j=1,2,\ldots,m$,%
\[
\nu\left(  d_{n,j}^{1}\right)  =n\mu_{1}+j\mu_{2}=\nu\left(  d_{n,j}%
^{2}\right)  .
\]
Thus, it follows from (\ref{eq-60}) that $\mathbf{B}^{\left(  \bm
d_{1}\right)  }=\mathbf{B}^{\left(  \bm    d_{2}\right)  }$, this also gives
that $\bm    \pi^{\left(  \bm    d_{1}\right)  }=\bm    \pi^{\left(  \bm
d_{2}\right)  }$. This completes the proof. \textbf{{\rule{0.08in}{0.08in}}}

Based on the problem description in Section 2, we define a suitable reward
function for the energy-efficient data center. For convenience of readers,
here we explain both energy consumption cost and system operational cost with
respect to some key factors of the data center. We summarize five classes of
costs in the data center as follows.

(a) \textit{Energy consumption cost}. It is seen from Section 2 that $P_{1,W}$
and $P_{2,W}$ are the energy consumption rates for the work state in Group 1
and Group 2, respectively; while $P_{2,S}$ is the energy consumption rate for
the sleep state only in Group 2. In addition, $C_{1}$ is the energy
consumption price.

(b) \textit{Holding cost}. Each job in the data center has to pay a holding
cost $C_{2}^{\left(  1\right)  }$ (resp. $C_{2}^{\left(  2\right)  }$) per
unit of sojourn time in Group 1 (resp. Group 2). Moreover, we have two assumed
conditions: A fast condition $\mu_{1}\geq\mu_{2},$ and a cheap condition
$C_{2}^{\left(  1\right)  }\leq C_{2}^{\left(  2\right)  }$.

(c) \textit{Transfer cost}. If there are some idle servers in Group 1, then
the jobs in the servers in Group 2 must be transferred to the idle servers in
Group 1 as many as possible. In this case, such each transfer needs to pay a
transfer cost $C_{3}$ for each job.

(d) \textit{Opportunity cost}. Once the data center has contained $n+m$ jobs,
then any arriving job has to be lost immediately. This leads to an opportunity
cost due to the job loss, hence $C_{4}$ is an opportunity cost for each lost job.

(e) \textit{Service price}. If a job finishes its service at a server and
leaves this system, then the data center gains a fixed service revenue (or
earnings) $R$ for each served job, that is, $R$ is the service price.

Based on the above cost and price definition, a reward function with respect
to both states and policies is defined as a profit rate (i.e. the total
revenues minus the total costs per unit of time). Therefore, the reward
function at state $\left(  I\left(  t\right)  ,J\left(  t\right)  \right)
^{\left(  \bm      d\right)  }$ under policy $\bm      d$ is defined as
\begin{align}
f^{\left(  \bm  d\right)  }\left(  i,j\right)  =  &  R\left[  i\mu_{1}+\left(
j\wedge d_{i,j}\right)  \mu_{2}\right]  -\left[  nP_{1,W}+d_{i,j}%
P_{2,W}+\left(  m-d_{i,j}\right)  P_{2,S}\right]  C_{1}\nonumber\\
&  -\left[  iC_{2}^{\left(  1\right)  }+jC_{2}^{\left(  2\right)  }\right]
-i\mu_{1}1_{\{j>0\}}C_{3}-\lambda1_{\left\{  i=n,j=m\right\}  }C_{4},
\label{eq-2}%
\end{align}
where $1_{\{\cdot\}}$ is an indicator function whose value is 1 when the event
in $\{\cdot\}$ happens, otherwise its value is 0. Furthermore, the job
transfer rate from Group 2 to Group 1 is given by $i\mu_{1}1_{\{j>0\}}$. If
$0\leq i\leq n-1$, then $j=0$ and $i\mu_{1}1_{\{j>0\}}C_{3}=0$. If $i=n$ and
$j=0$, then $i\mu_{1}1_{\{j>0\}}C_{3}=0$. If $i=n$ and $1\leq j\leq m$, then
$i\mu_{1}1_{\{j>0\}}C_{3}=n\mu_{1}C_{3}$.

For convenience of readers, it is necessary to explain the reward function
from four different cases as follows.

\textbf{Case (a):} For $i=0$ and $j=0$,
\begin{equation}
f\left(  0,0\right)  =-\left(  nP_{1,W}+mP_{2,S}\right)  C_{1}. \label{eq-50}%
\end{equation}

\textbf{Case (b): }For $i=1,2,\ldots,n$ and $j=0$,
\begin{equation}
f\left(  i,0\right)  =Ri\mu_{1}-\left(  nP_{1,W}+mP_{2,S}\right)  C_{1}%
-iC_{2}^{\left(  1\right)  }. \label{eq-49}%
\end{equation}

Note that in Cases (a) and (b), there is no job in Group 2, thus each server
of Group 2 is at the sleep state. However, in the following two cases (c) and
(d), there are some jobs in Group 2, hence the policy $\bm      d$\ will play
a key role in opening or closing some servers of Group 2 in order to
effectively save energy.

\textbf{Case (c): }For $i=n$, $j=1,2,\ldots,m-1$; and $d_{n,j}=0,1,\ldots,m$,
\begin{align}
f^{\left(  \bm  d\right)  }\left(  n,j\right)   &  =R\left[  n\mu_{1}+\left(
j\wedge d_{n,j}\right)  \mu_{2}\right]  -\left[  nP_{1,W}+d_{n,j}%
P_{2,W}+\left(  m-d_{n,j}\right)  P_{2,S}\right]  C_{1}\nonumber\\
&  \text{ \ \ }-\left[  nC_{2}^{\left(  1\right)  }+jC_{2}^{\left(  2\right)
}\right]  -n\mu_{1}C_{3}. \label{eq-51}%
\end{align}
To further simplify or compute (\ref{eq-51}), we need to especially deal with
$j\wedge d_{n,j}$. To this end, if $d_{n,j}=1,2,\ldots,j$, then $j\wedge
d_{n,j}=d_{n,j}$, hence we have%
\begin{align}
f^{\left(  \bm  d\right)  }\left(  n,j\right)   &  =Rn\mu_{1}+\left[  R\mu
_{2}-\left(  P_{2,W}-P_{2,S}\right)  C_{1}\right]  d_{n,j}\nonumber\\
&  \text{ \ \ }-\left(  nP_{1,W}+mP_{2,S}\right)  C_{1}-\left[  nC_{2}%
^{\left(  1\right)  }+jC_{2}^{\left(  2\right)  }\right]  -n\mu_{1}C_{3};
\label{eq-52}%
\end{align}
while if $d_{n,j}=j+1,j+2,\ldots,m$, then $j\wedge d_{n,j}=j$, hence we have%
\begin{align}
f^{\left(  \bm  d\right)  }\left(  n,j\right)   &  =R\left(  n\mu_{1}+j\mu
_{2}\right)  -\left(  P_{2,W}-P_{2,S}\right)  C_{1}d_{n,j}\nonumber\\
&  \text{ \ \ }-\left(  nP_{1,W}+mP_{2,S}\right)  C_{1}-\left[  nC_{2}%
^{\left(  1\right)  }+jC_{2}^{\left(  2\right)  }\right]  -n\mu_{1}C_{3}.
\label{eq-53}%
\end{align}

\textbf{Case (d): }For $i=n$ and $j=m$; and $d_{n,m}=0,1,\ldots,m$, we obtain
that $m\wedge d_{n,m}=d_{n,m}$, and we have
\begin{align}
f^{\left(  \bm  d\right)  }\left(  n,m\right)   &  =R\left(  n\mu_{1}%
+d_{n,m}\mu_{2}\right)  -\left[  nP_{1,W}+d_{n,m}P_{2,W}+\left(
m-d_{n,m}\right)  P_{2,S}\right]  C_{1}\nonumber\\
&  \text{ \ \ }-\left[  nC_{2}^{\left(  1\right)  }+mC_{2}^{\left(  2\right)
}\right]  -n\mu_{1}C_{3}-\lambda C_{4}\nonumber\\
&  =Rn\mu_{1}+\left[  R\mu_{2}-\left(  P_{2,W}-P_{2,S}\right)  C_{1}\right]
d_{n,m}\nonumber\\
&  \text{ \ \ }-\left(  nP_{1,W}+mP_{2,S}\right)  C_{1}-\left[  nC_{2}%
^{\left(  1\right)  }+mC_{2}^{\left(  2\right)  }\right]  -n\mu_{1}%
C_{3}-\lambda C_{4}. \label{eq-54}%
\end{align}

We define an $\left(  n+m+1\right)  $-dimensional column vector composed of
the elements $f\left(  i,j\right)  $ and $f^{\left(  \bm      d\right)
}\left(  i,j\right)  $\ as
\begin{equation}
\bm  f^{\left(  \bm  d\right)  }=\left(  f\left(  0,0\right)  ,f\left(
1,0\right)  ,\ldots,f\left(  n,0\right)  ,f^{\left(  \bm  d\right)  }\left(
n,1\right)  ,\ldots,f^{\left(  \bm  d\right)  }\left(  n,m\right)  \right)
^{T}, \label{eq-33}%
\end{equation}
where $a^{T}$ denotes transpose of vector or matrix $a$.

In the remainder of this section, the long-run average profit of the data
center (or the policy-based continuous-time birth-death process $\left\{
\mathbf{X}^{\left(  \bm      d\right)  }(t):t\geq0\right\}  $) under an
energy-efficient policy $\bm      d$ is defined as%
\begin{align}
\eta^{\bm  d}  &  =\lim_{T\rightarrow+\infty}E\left\{  \frac{1}{T}\int_{0}%
^{T}f^{(\bm  d)}\left(  \left(  I\left(  t\right)  ,J\left(  t\right)
\right)  ^{(\bm  d)}\right)  \text{d}t\right\} \nonumber\\
&  =\lim_{T\rightarrow+\infty}E\left\{  \frac{1}{T}\int_{0}^{T}f^{(\bm
d)}\left(  \mathbf{X}^{\left(  \bm  d\right)  }(t)\right)  \text{d}t\right\}
\nonumber\\
&  =\bm  \pi^{\left(  \bm  d\right)  }\bm  f^{\left(  \bm  d\right)  },
\label{eq-6}%
\end{align}
where $\bm      \pi^{\left(  \bm      d\right)  }$ and $\bm      f^{\left(
\bm  d\right)  }$ are given by (\ref{eq-1-1}) and (\ref{eq-33}), respectively.

We observe that as the number of working servers in Group 2 decreases, the
total revenues and the total costs in the data center will decrease
synchronously, vice versa. Thus, there is a tradeoff between the total
revenues and the total costs. This motivates us to study an optimal mechanism
design for the energy-efficient data center. The objective is to find an
optimal energy-efficient policy $\bm      d^{\ast}$ such that the long-run
average profit $\eta^{\bm      d}$ is maximize, that is,%
\begin{equation}
\bm  d^{\ast}=\underset{\bm  d\in\mathcal{D}}{\arg\max}\left\{  \eta^{\bm
d}\right\}  . \label{eq-7}%
\end{equation}

In fact, it is difficult and challenging to analyze the properties of the
optimal energy-efficient policy $\bm      d^{\ast}$, and to provide an
effective algorithm for computing the optimal policy $\bm      d^{\ast}$. In
the next section, we will introduce the sensitivity-based optimization theory
to study this energy-efficient optimization problem.

\section{The Poisson Equation and Its Explicit Solution}

In this section, for the energy-efficient data center, we set up a Poisson
equation which is derived by means of the law of total probability focusing on
some stop times. It is worth noting that the Poisson equation provides a
useful relation between the sensitivity-based optimization and the Markov
decision processes (MDPs). Also, we use the RG-factorization, given in Li and
Cao \cite{Li:2004} or Li \cite{Li:2010}, to solve the Poisson equation and
provide the explicit expression for its unique solution.

For $\bm     d\in\mathcal{D}$, it follows from Chapter 2 in Cao
\cite{Cao:2007} that for the policy-based continuous-time Markov process
$\left\{  \mathbf{X}^{\left(  \bm     d\right)  }(t):t\geq0\right\}  $, we
define the performance potential as%
\begin{align}
g^{\left(  \bm  d\right)  }\left(  i,j\right)   &  =E\left\{  \left.  \int
_{0}^{+\infty}\left[  f^{(\bm  d)}\left(  \left(  I\left(  t\right)  ,J\left(
t\right)  \right)  ^{\left(  \bm  d\right)  }\right)  -\eta^{\bm  d}\right]
\text{d}t\right|  \left(  I\left(  0\right)  ,J\left(  0\right)  \right)
^{\left(  \bm  d\right)  }=\left(  i,j\right)  \right\} \nonumber\\
&  =E\left\{  \left.  \int_{0}^{+\infty}\left[  f^{(\bm  d)}\left(
\mathbf{X}^{\left(  \bm  d\right)  }(t)\right)  -\eta^{\bm  d}\right]
\text{d}t\right|  \mathbf{X}^{\left(  \bm  d\right)  }\left(  0\right)
=\left(  i,j\right)  \right\}  , \label{eq-8}%
\end{align}
where $\eta^{\bm     d}$ is defined in (\ref{eq-6}). It is seen from Cao
\cite{Cao:2007} that for any policy $\bm     d\in\mathcal{D}$, $g^{\left(
\bm  d\right)  }\left(  i,j\right)  $ quantifies the contribution of the
initial state $(i,j)$ to the long-run average profit of the data center. Here
$g^{\left(  \bm     d\right)  }\left(  i,j\right)  $ is also called the
relative value function or the bias in the traditional MDP theory, see, e.g.
Puterman \cite{Put:2014}. We further define a column vector $\bm
g^{\left(  \bm   d\right)  }$ with elements $g^{\left(  \bm     d\right)
}\left(  i,j\right)  $ for $(i,j)\in\bm     \Omega$ as
\begin{equation}
\bm  g^{\left(  \bm  d\right)  }=\left(  g^{\left(  \bm  d\right)  }\left(
0,0\right)  ,g^{\left(  \bm  d\right)  }\left(  1,0\right)  ,\ldots,g^{\left(
\bm  d\right)  }\left(  n,0\right)  ,g^{\left(  \bm  d\right)  }\left(
n,1\right)  ,\ldots,g^{\left(  \bm  d\right)  }\left(  n,m\right)  \right)
^{T}. \label{eq-34}%
\end{equation}

We define the first departure time from state $\left(  i,j\right)  $ as
\[
\tau=\inf\left\{  t\geq0:\left(  I\left(  t\right)  ,J\left(  t\right)
\right)  ^{\left(  \bm d\right)  }\neq\left(  i,j\right)  \right\}  ,
\]
where $\left(  I\left(  0\right)  ,J\left(  0\right)  \right)  ^{\left(  \bm
d\right)  }=\left(  i,j\right)  $. Clearly, $\tau$ is a stop time of the
Markov process $\left\{  \mathbf{X}^{\left(  \bm d\right)  }(t):t\geq
0\right\}  $. Based on this, if $\left(  i,j\right)  =\left(  0,0\right)  $,
then it is seen from\ (\ref{eq-60}) that state $\left(  0,0\right)  $ is the
upper boundary state of the birth-death process $\mathbf{B}^{\left(  \bm
d\right)  }$, hence $\left(  I\left(  \tau\right)  ,J\left(  \tau\right)
\right)  ^{\left(  \bm d\right)  }=\left(  1,0\right)  $. Similarly, we get a
basic relation for the stop time $\tau$ as follows.%
\begin{equation}
\left(  I\left(  \tau\right)  ,J\left(  \tau\right)  \right)  ^{\left(  \bm
d\right)  }=\left\{
\begin{array}
[c]{ll}%
\left(  1,0\right)  , & \left(  i,j\right)  =\left(  0,0\right)  ,\\
\left(  i-1,0\right)  ,\left(  i+1,0\right)  , & \left(  i,j\right)  =\left(
i,0\right)  ,i=1,2,\ldots,n-1,\\
\left(  n-1,0\right)  ,\left(  n,1\right)  , & \left(  i,j\right)  =\left(
n,0\right)  ,\\
\left(  n,j-1\right)  ,\left(  n,j+1\right)  , & \left(  i,j\right)  =\left(
n,j\right)  ,j=1,2,\ldots,m-1,\\
\left(  n,m-1\right)  , & \left(  i,j\right)  =\left(  n,m\right)  .
\end{array}
\right.  \label{eq-ij}%
\end{equation}

Now, we derive a Poisson equation to compute the column vector $\bm
g^{\left(  \bm                                   d\right)  }$ in terms of the
stop time $\tau$ and the basic relation (\ref{eq-ij}). By a similar
computation to that in Li and Cao \cite{Li:2004} or Xia et al.
\cite{Xia:2017a}, our analysis is decomposed into five parts as follows.

For $i=1,2,\ldots,n-1\ $and $j=0,$ we have%
\begin{align}
g^{\left(  \bm d\right)  }\left(  i,0\right)  = &  E\left\{  \left.  \int
_{0}^{+\infty}\left[  f^{\left(  \bm d\right)  }\left(  \mathbf{X}^{\left(
\bm d\right)  }(t)\right)  -\eta^{\bm d}\right]  \text{d}t\right|  \left(
I\left(  0\right)  ,J\left(  0\right)  \right)  ^{\left(  \bm d\right)
}=\left(  i,0\right)  \right\}  \nonumber\\
= &  E\left\{  \left.  \tau\right|  I\left(  t\right)  =i,J\left(  t\right)
=0\right\}  \left[  f\left(  i,0\right)  -\eta^{\bm d}\right]  \nonumber\\
&  +E\left\{  \left.  \int_{\tau}^{+\infty}\left[  f^{\left(  \bm d\right)
}\left(  \mathbf{X}^{\left(  \bm d\right)  }(t)\right)  -\eta^{\bm d}\right]
\text{d}t\right|  \left(  I\left(  \tau\right)  ,J\left(  \tau\right)
\right)  ^{\left(  \bm d\right)  }\right\}  \nonumber\\
= &  \frac{1}{\lambda+i\mu_{1}}\left[  f\left(  i,0\right)  -\eta^{\bm
d}\right]  \nonumber\\
&  +\frac{\lambda}{\lambda+i\mu_{1}}E\left\{  \left.  \int_{0}^{+\infty
}\left[  f^{\left(  \bm d\right)  }\left(  \mathbf{X}^{\left(  \bm d\right)
}(t)\right)  -\eta^{\bm d}\right]  \text{d}t\right|  \left(  I\left(
0\right)  ,J\left(  0\right)  \right)  ^{\left(  \bm d\right)  }=\left(
i+1,0\right)  \right\}  \nonumber\\
&  +\frac{i\mu_{1}}{\lambda+i\mu_{1}}E\left\{  \left.  \int_{0}^{+\infty
}\left[  f^{\left(  \bm d\right)  }\left(  \mathbf{X}^{\left(  \bm d\right)
}(t)\right)  -\eta^{\bm d}\right]  \text{d}t\right|  \left(  I\left(
0\right)  ,J\left(  0\right)  \right)  ^{\left(  \bm d\right)  }=\left(
i-1,0\right)  \right\}  \nonumber\\
= &  \frac{1}{\lambda+i\mu_{1}}\left[  f\left(  i,0\right)  -\eta^{\bm
d}\right]  +\frac{\lambda}{\lambda+i\mu_{1}}g^{\left(  \bm d\right)  }\left(
i+1,0\right)  +\frac{i\mu_{1}}{\lambda+i\mu_{1}}g^{\left(  \bm d\right)
}\left(  i-1,0\right)  ,\label{eq-46}%
\end{align}
where for the birth-death process $\left\{  \mathbf{X}^{\left(  \bm d\right)
}(t):t\geq0\right\}  $, it is easy to see that%
\[
\int_{0}^{\tau}\left[  f^{\left(  \bm d\right)  }\left(  \left(  I\left(
t\right)  ,J\left(  t\right)  \right)  ^{(\bm d)}\right)  -\eta^{\bm
d}\right]  \text{d}t=\tau\left[  f\left(  i,0\right)  -\eta^{\bm d}\right]  ,
\]%
\[
E\left\{  \left.  \tau\right|  I\left(  t\right)  =i,J\left(  t\right)
=0\right\}  =\frac{1}{\lambda+i\mu_{1}}.
\]
Thus, we obtain%
\begin{equation}
i\mu_{1}g^{\left(  \bm d\right)  }\left(  i-1,0\right)  -\left(  \lambda
+i\mu_{1}\right)  g^{\left(  \bm d\right)  }\left(  i,0\right)  +\lambda
g^{\left(  \bm d\right)  }\left(  i+1,0\right)  =\eta^{\bm d}-f\left(
i,0\right)  .\label{eq-46-1}%
\end{equation}
Base on (\ref{eq-46}), with a boundary consideration, for $i=0\ $and $j=0$, we
have
\begin{equation}
-\lambda g^{\left(  \bm d\right)  }\left(  0,0\right)  +\lambda g^{\left(  \bm
d\right)  }\left(  1,0\right)  =\eta^{\bm d}-f\left(  0,0\right)
.\label{eq-45}%
\end{equation}
For $i=n\ $and $j=0,$ we have%
\begin{equation}
n\mu_{1}g^{\left(  \bm d\right)  }\left(  n-1,0\right)  -\left(  \lambda
+n\mu_{1}\right)  g^{\left(  \bm d\right)  }\left(  n,0\right)  +\lambda
g^{\left(  \bm d\right)  }\left(  n,1\right)  =\eta^{\bm d}-f\left(
n,0\right)  .\label{eq-47-1}%
\end{equation}

For $i=n\ $and $j=1,2,\ldots,m-1$, we have%
\begin{align}
g^{\left(  \bm  d\right)  }\left(  n,j\right)  =  &  E\left\{  \left.
\int_{0}^{+\infty}\left[  f^{\left(  \bm  d\right)  }\left(  \mathbf{X}%
^{\left(  \bm  d\right)  }(t)\right)  -\eta^{\bm  d}\right]  \text{d}t\right|
\left(  I\left(  0\right)  ,J\left(  0\right)  \right)  ^{\left(  \bm
d\right)  }=\left(  n,j\right)  \right\} \nonumber\\
=  &  E\left\{  \left.  \tau\right|  I\left(  t\right)  =n,J\left(  t\right)
=j\right\}  \left[  f^{\left(  \bm  d\right)  }\left(  n,j\right)  -\eta^{\bm
d}\right] \nonumber\\
&  +E\left\{  \left.  \int_{\tau}^{+\infty}\left[  f^{\left(  \bm  d\right)
}\left(  \mathbf{X}^{\left(  \bm  d\right)  }(t)\right)  -\eta^{\bm
d}\right]  \text{d}t\right|  \left(  I\left(  \tau\right)  ,J\left(
\tau\right)  \right)  ^{\left(  \bm  d\right)  }\right\} \nonumber\\
=  &  \frac{1}{\lambda+\nu\left(  d_{n,j}\right)  }\left[  f^{\left(  \bm
d\right)  }\left(  n,j\right)  -\eta^{\bm  d}\right] \nonumber\\
&  +\frac{\lambda}{\lambda+\nu\left(  d_{n,j}\right)  }E\left\{  \left.
\int_{0}^{+\infty}\left[  f^{\left(  \bm  d\right)  }\left(  \mathbf{X}%
^{\left(  \bm  d\right)  }(t)\right)  -\eta^{\bm  d}\right]  \text{d}t\right|
\left(  I\left(  0\right)  ,J\left(  0\right)  \right)  ^{\left(  \bm
d\right)  }=\left(  n,j+1\right)  \right\} \nonumber\\
&  +\frac{\nu\left(  d_{n,j}\right)  }{\lambda+\nu\left(  d_{n,j}\right)
}E\left\{  \left.  \int_{0}^{+\infty}\left[  f^{\left(  \bm  d\right)
}\left(  \mathbf{X}^{\left(  \bm  d\right)  }(t)\right)  -\eta^{\bm
d}\right]  \text{d}t\right|  \left(  I\left(  0\right)  ,J\left(  0\right)
\right)  ^{\left(  \bm  d\right)  }=\left(  n,j-1\right)  \right\} \nonumber\\
=  &  \frac{1}{\lambda+\nu\left(  d_{n,j}\right)  }\left[  f^{\left(  \bm
d\right)  }\left(  n,j\right)  -\eta^{\bm  d}\right]  +\frac{\lambda}%
{\lambda+\nu\left(  d_{n,j}\right)  }g^{\left(  \bm  d\right)  }\left(
n,j+1\right) \nonumber\\
&  +\frac{\nu\left(  d_{n,j}\right)  }{\lambda+\nu\left(  d_{n,j}\right)
}g^{\left(  \bm  d\right)  }\left(  n,j-1\right)  , \label{eq-10}%
\end{align}
where%
\begin{equation}
E\left\{  \left.  \tau\right|  I\left(  t\right)  =n,J\left(  t\right)
=j\right\}  =\frac{1}{\lambda+\nu\left(  d_{n,j}\right)  }. \label{eq-11}%
\end{equation}
It follows from (\ref{eq-10}) that%
\begin{equation}
\nu\left(  d_{n,j}\right)  g^{\left(  \bm  d\right)  }\left(  n,j-1\right)
-\left[  \lambda+\nu\left(  d_{n,j}\right)  \right]  g^{\left(  \bm  d\right)
}\left(  n,j\right)  +\lambda g^{\left(  \bm  d\right)  }\left(  n,j+1\right)
=\eta^{\bm  d}-f^{\left(  \bm  d\right)  }\left(  n,j\right)  .
\label{eq-10-1}%
\end{equation}
For $i=n\ $and $j=m$, with a boundary consideration, a similar analysis to
(\ref{eq-10-1}) gives%
\begin{equation}
\nu\left(  d_{n,m}\right)  g^{\left(  \bm  d\right)  }\left(  n,m-1\right)
-\nu\left(  d_{n,m}\right)  g^{\left(  \bm  d\right)  }\left(  n,m\right)
=\eta^{\bm  d}-f^{\left(  \bm  d\right)  }\left(  n,m\right)  .
\label{eq-48-1}%
\end{equation}
Note that $\nu\left(  d_{n,m}\right)  =n\mu_{1}+\left(  d_{n,m}\wedge
m\right)  \mu_{2}=n\mu_{1}+d_{n,m}\mu_{2}$ due to the fact that $d_{n,m}%
\in\left\{  0,1,\ldots,m\right\}  $.

It follows from (\ref{eq-46-1}), (\ref{eq-45}), (\ref{eq-47-1}),
(\ref{eq-10-1}) and (\ref{eq-48-1}) that
\[
\mathbf{B}^{\left(  \bm                                  d\right)  }\bm
g^{(\bm     d)}=\eta^{\bm                             d}\bm      e-\bm
f^{(\bm            d)}%
\]
or%
\begin{equation}
-\mathbf{B}^{\left(  \bm  d\right)  }\bm  g^{(\bm  d)}=\bm  f^{\left(  \bm
d\right)  }-\eta^{\bm  d}\bm  e, \label{eq-61}%
\end{equation}
where $\bm                                  f^{\left(  \bm
d\right)  }$ is given in (\ref{eq-33}) and $\mathbf{B}^{\left(  \bm
d\right)  }$ is given in (\ref{eq-60}).

To solve the system of linear equations (\ref{eq-61}), we note that
rank$\left(  \mathbf{B}^{\left(  \bm    d\right)  }\right)  =n+m$ and
$\det\left(  \mathbf{B}^{\left(  \bm    d\right)  }\right)  =0$ due to that
the size of the matrix $\mathbf{B}^{\left(  \bm    d\right)  }$ is $n+m+1$.
Hence, this system of linear equations (\ref{eq-61}) exists infinitely-many
solutions with a free constant of an additive term. Let $\mathcal{B}$ be a
matrix obtained through omitting the first row and the first column vectors of
the matrix $\mathbf{B}^{\left(  \bm    d\right)  }$. Then,%
\[
\mathcal{B}=\left(
\begin{array}
[c]{cccccc}%
-\left(  \lambda+\mu_{1}\right)  & \lambda &  &  &  & \\
2\mu_{1} & -\left(  \lambda+2\mu_{1}\right)  & \lambda &  &  & \\
& \ddots\text{ \ \ \ \ \ \ } & \ddots\text{ \ \ \ \ \ \ \ } & \ddots\text{
\ \ \ \ \ \ } &  & \\
& n\mu_{1} & -\left(  \lambda+n\mu_{1}\right)  & \lambda &  & \\
&  & \nu\left(  d_{n,1}\right)  & -\left[  \lambda+\nu\left(  d_{n,1}\right)
\right]  & \lambda & \\
&  & \text{\ \ \ \ \ \ \ }\ddots & \ \ \ \ \ \ \  & \ddots\text{
\ \ \ \ \ \ \ \ \ \ }\ddots & \text{ \ \ \ \ \ }\\
&  &  & \nu\left(  d_{n,m-1}\right)  & -\left[  \lambda+\nu\left(
d_{n,m-1}\right)  \right]  & \lambda\\
&  &  &  & \nu\left(  d_{n,m}\right)  & -\nu\left(  d_{n,m}\right)
\end{array}
\right)  .
\]
Obviously, rank$\left(  \mathcal{B}\right)  =n+m$ and the size of the matrix
$\mathcal{B}$ is $n+m$. Hence, the matrix $\mathcal{B}$ is invertible.

Let $\bm{h} ^{\left(  \bm d\right)  }$ and $\bm \varphi^{(\bm d)}$ be two
column vectors of size $n+m$ obtained through omitting the first element of
the two column vectors $\bm f^{\left(  \bm d\right)  }-\eta^{\bm d}\bm e$ and
$\bm g^{(\bm d)}$ with size $n+m+1$, respectively. Then,%
\[
\bm{h} ^{\left(  \bm d\right)  }=\left(
\begin{array}
[c]{c}%
f\left(  1,0\right)  -\eta^{\bm d}\\
\vdots\\
f\left(  n,0\right)  -\eta^{\bm d}\\
f^{\left(  \bm d\right)  }\left(  n,1\right)  -\eta^{\bm d}\\
\vdots\\
f^{\left(  \bm d\right)  }\left(  n,m\right)  -\eta^{\bm d}%
\end{array}
\right)  \overset{\text{def}}{=}\left(
\begin{array}
[c]{c}%
h_{1}^{\left(  \bm d\right)  }\\
\vdots\\
h_{n}^{\left(  \bm d\right)  }\\
h_{n+1}^{\left(  \bm d\right)  }\\
\vdots\\
h_{n+m}^{\left(  \bm d\right)  }%
\end{array}
\right)  ,\text{ \ }\bm \varphi^{(\bm d)}=\left(
\begin{array}
[c]{c}%
g^{\left(  \bm d\right)  }\left(  1,0\right)  \\
\vdots\\
g^{\left(  \bm d\right)  }\left(  n,0\right)  \\
g^{\left(  \bm d\right)  }\left(  n,1\right)  \\
\vdots\\
g^{\left(  \bm d\right)  }\left(  n,m\right)
\end{array}
\right)  .
\]
Therefore, it follows from (\ref{eq-61}) that%
\begin{equation}
-\mathcal{B}\bm \varphi^{(\bm d)}=\bm{h} ^{\left(  \bm d\right)  }+\mu
_{1}\bm{e} _{1}g^{\left(  \bm d\right)  }\left(  0,0\right)  ,\label{Eq-1}%
\end{equation}
where $\bm{e} _{1}$ is a column vector with the first element be one and all
the others be zero. Note that the matrix $-\mathcal{B}$ is invertible and
$\left(  -\mathcal{B}\right)  ^{-1}>0$, thus, the system of linear equations
(\ref{Eq-1}) always exists one unique solution%
\begin{equation}
\bm \varphi^{(\bm d)}=\left(  -\mathcal{B}\right)  ^{-1}\bm{h} ^{\left(  \bm
d\right)  }+\mu_{1}\left(  -\mathcal{B}\right)  ^{-1}\bm{e} _{1}\cdot
\Im,\label{Eq-2}%
\end{equation}
where $g^{\left(  \bm d\right)  }\left(  0,0\right)  =\Im$ is any given
constant. For convenience of computation, we take $g^{\left(  \bm d\right)
}\left(  0,0\right)  =\Im=1$. In this case, we have%
\begin{equation}
\bm \varphi^{(\bm d)}=\left(  -\mathcal{B}\right)  ^{-1}\bm{h} ^{\left(  \bm
d\right)  }+\mu_{1}\left(  -\mathcal{B}\right)  ^{-1}\bm{e} _{1}%
.\label{Eq-2-1}%
\end{equation}

To provide an explicit expression for the unique solution to the system of
linear equations (\ref{Eq-1}), it is easy to see from (\ref{Eq-2-1}) that we
need to first establish an explicit expression for the invertible matrix
$\left(  -\mathcal{B}\right)  ^{-1}$. While the explicit expression of the
invertible matrix $\left(  -\mathcal{B}\right)  ^{-1}$ can be obtained by
means of the RG-factorization, which is given in Li and Cao \cite{Li:2004} for
QBD processes, and more generally, Li \cite{Li:2010} for general Markov processes.

To express the invertible matrix $\left(  -\mathcal{B}\right)  ^{-1}$, we
write the UL-type $U$-measure as%
\[%
\begin{array}
[c]{ll}
& U_{n+m}=-\nu\left(  d_{n,m}\right)  ,\\
& U_{n+m-1}=-\left[  \lambda+\nu\left(  d_{n,m-1}\right)  \right]
+\lambda\left(  -U_{n+m}\right)  ^{-1}\nu\left(  d_{n,m}\right)  ,\\
& U_{k}=A_{1}^{\left(  k\right)  }+\lambda\left(  -U_{k+1}\right)  ^{-1}%
A_{2}^{\left(  k+1\right)  },\text{ \ \ }k=n+m-2,n+m-3,\ldots,1,
\end{array}
\]
and the UL-type $R$- and $G$-measures as
\[%
\begin{array}
[c]{lll}%
R_{k} & = & \lambda\left(  -U_{k+1}\right)  ^{-1},\text{ \ }k=1,2,\ldots
,n+m-1,\\
G_{k} & = & \left(  -U_{k}\right)  ^{-1}A_{2}^{\left(  k\right)  },\text{
\ }k=2,3,\ldots,n+m,
\end{array}
\]
where
\[%
\begin{array}
[c]{ll}
& A_{1}^{\left(  k\right)  }=\left\{
\begin{array}
[c]{ll}%
-\left(  \lambda+k\mu_{1}\right)  , & k=1,2,\ldots,n,\\
-\left[  \lambda+\nu\left(  d_{n,k-n}\right)  \right]  , & k=n+1,n+2,\ldots
,n+m-1,
\end{array}
\right. \\
& A_{2}^{\left(  k\right)  }=\left\{
\begin{array}
[c]{ll}%
k\mu_{1}, & k=2,3,\ldots,n,\\
\nu\left(  d_{n,k-n}\right)  , & k=n+1,n+2,\ldots,n+m.
\end{array}
\right.
\end{array}
\]
Thus, the UL-type RG-factorization of the birth-death process $\mathcal{B}$ is
given by
\[
\mathcal{B}=\left(  I-R_{U}\right)  U_{D}\left(  I-G_{L}\right)  ,
\]
where%
\[
R_{U}=\left(
\begin{array}
[c]{ccccc}%
0 & R_{1} &  &  & \\
& 0 & R_{2} &  & \\
&  & \ddots & \ddots & \\
&  &  & 0 & R_{n+m-1}\\
&  &  &  & 0
\end{array}
\right)  ,\text{ \ }G_{L}=\left(
\begin{array}
[c]{ccccc}%
0 &  &  &  & \\
G_{2} & 0 &  &  & \\
& G_{3} & 0 &  & \\
&  & \ddots & \ddots & \\
&  &  & G_{n+m} & 0
\end{array}
\right)  ,
\]
and%
\[
U_{D}=\text{diag}\left(  U_{1},U_{2},\ldots,U_{n+m}\right)  .
\]
Therefore, we obtain%
\[
\left(  -\mathcal{B}\right)  ^{-1}=\left(  I-G_{L}\right)  ^{-1}\left(
-U_{D}\right)  ^{-1}\left(  I-R_{U}\right)  ^{-1}.
\]

Let%
\begin{align*}
X_{k}^{\left(  l\right)  }  &  =R_{l}R_{l+1}\cdots R_{l+k-2},\text{ }1\leq
l\leq n+m-1,2\leq k\leq n+m,\\
Y_{k}^{\left(  l\right)  }  &  =G_{l}G_{l-1}\cdots G_{l-k+2},\text{ }2\leq
k\leq l\leq n+m.
\end{align*}
Then,%
\[
\left(  I-G_{L}\right)  ^{-1}=\left(
\begin{array}
[c]{ccccc}%
1 &  &  &  & \\
Y_{2}^{\left(  2\right)  } & 1 &  &  & \\
Y_{3}^{\left(  3\right)  } & Y_{2}^{\left(  3\right)  } & 1 &  & \\
\vdots & \vdots & \ddots & \ddots & \\
Y_{n+m}^{\left(  n+m\right)  } & Y_{n+m-1}^{\left(  n+m\right)  } & \cdots &
Y_{2}^{\left(  n+m\right)  } & 1
\end{array}
\right)  ,
\]%
\[
\left(  I-R_{U}\right)  ^{-1}=\left(
\begin{array}
[c]{ccccc}%
1 & X_{2}^{\left(  1\right)  } & X_{3}^{\left(  1\right)  } & \cdots &
X_{n+m}^{\left(  1\right)  }\\
& 1 & X_{2}^{\left(  2\right)  } & \cdots & X_{n+m-1}^{\left(  2\right)  }\\
&  & 1 & \cdots & X_{n+m-2}^{\left(  3\right)  }\\
&  &  & \ddots & \vdots\\
&  &  &  & 1
\end{array}
\right)  ,
\]
and%
\[
\left(  -U_{D}\right)  ^{-1}=\text{diag}\left(  \left(  -U_{1}\right)
^{-1},\left(  -U_{2}\right)  ^{-1},\ldots,\left(  -U_{n+m}\right)
^{-1}\right)  .
\]
Thus, we obtain the explicit expression%
\begin{align}
\left(  -\mathcal{B}\right)  ^{-1}  &  =\left(
\begin{array}
[c]{ccccc}%
1 &  &  &  & \\
Y_{2}^{\left(  2\right)  } & 1 &  &  & \\
Y_{3}^{\left(  3\right)  } & Y_{2}^{\left(  3\right)  } & 1 &  & \\
\vdots & \vdots & \ddots & \ddots & \\
Y_{n+m}^{\left(  n+m\right)  } & Y_{n+m-1}^{\left(  n+m\right)  } & \cdots &
Y_{2}^{\left(  n+m\right)  } & 1
\end{array}
\right)  \text{diag}\left(  \left(  -U_{1}\right)  ^{-1},\left(
-U_{2}\right)  ^{-1},\ldots,\left(  -U_{n+m}\right)  ^{-1}\right) \nonumber\\
&  \cdot\left(
\begin{array}
[c]{ccccc}%
1 & X_{2}^{\left(  1\right)  } & X_{3}^{\left(  1\right)  } & \cdots &
X_{n+m}^{\left(  1\right)  }\\
& 1 & X_{2}^{\left(  2\right)  } & \cdots & X_{n+m-1}^{\left(  2\right)  }\\
&  & 1 & \cdots & X_{n+m-2}^{\left(  3\right)  }\\
&  &  & \ddots & \vdots\\
&  &  &  & 1
\end{array}
\right)  . \label{Eq-3}%
\end{align}

It is worthwhile to note that for the $U$-, $R$- and $G$-measures discussed
above, each of them depends on the policy $\bm      d\in\mathcal{D}$. For
simplification of descriptions, we omit such superscript $\bm      d$ in their
formulae or equations.

\begin{Rem}
In general, it is difficult to solve the high-dimensional system of linear
equations (\ref{Eq-1}) for a general Markov process. However, the
RG-factorization method is usually very effective for solving such problem for
a level-dependent QBD process, in which the expression of the invertible
matrix $\left(  -\mathcal{B}\right)  ^{-1}$ is similar to that in (\ref{Eq-3})
with elements changed from scalar to block matrices. In addition, it is also a
key that the vector $\bm      \varphi^{(\bm      d)}$ can be numerically
computed from (\ref{Eq-2}) and (\ref{Eq-3}) by using the RG-factorization. See
Li and Cao \cite{Li:2004} for more details.
\end{Rem}

The following theorem provides an explicit expression for the vector $\bm
\varphi^{(\bm                         d)}$ under a constraint condition
$g^{\left(  \bm           d\right)  }\left(  0,0\right)  =1$. Note that this
expression is very useful for applications of the sensitivity-based
optimization theory in our later study.

\begin{The}
\label{The:Potential}If $g^{\left(  \bm   d\right)  }\left(  0,0\right)  =1$,
then for $k=1,$%
\[
g^{\left(  \bm   d\right)  }\left(  1,0\right)  =\left(  -U_{1}\right)
^{-1}\left[  h_{1}^{\left(  \bm   d\right)  }+\sum_{j=2}^{m+n}X_{j}^{\left(
1\right)  }h_{j}^{\left(  \bm   d\right)  }\right]  +G_{1};
\]
For $2\leq k\leq n$,%
\begin{align*}
g^{\left(  \bm  d\right)  }\left(  k,0\right)  =  &  \left(  -U_{k}\right)
^{-1}\left[  h_{k}^{\left(  \bm  d\right)  }+\sum_{j=2}^{m+n-k+1}%
X_{j}^{\left(  k\right)  }h_{j+k-1}^{\left(  \bm  d\right)  }\right] \\
&  +\sum_{i=1}^{k-1}Y_{k-i+1}^{\left(  k\right)  }\left(  -U_{i}\right)
^{-1}\left[  h_{i}^{\left(  \bm  d\right)  }+\sum_{j=2}^{m+n-i+1}%
X_{j}^{\left(  i\right)  }h_{j+i-1}^{\left(  \bm  d\right)  }\right]
+\prod_{j=1}^{k}G_{j};
\end{align*}
For $n+1\leq k\leq n+m-1$,%
\begin{align*}
g^{\left(  \bm  d\right)  }\left(  n,k-n\right)  =  &  \left(  -U_{k}\right)
^{-1}\left[  h_{k}^{\left(  \bm  d\right)  }+\sum_{j=2}^{m+n-k+1}%
X_{j}^{\left(  k\right)  }h_{j+k-1}^{\left(  \bm  d\right)  }\right] \\
&  +\sum_{i=1}^{k-1}Y_{k-i+1}^{\left(  k\right)  }\left(  -U_{i}\right)
^{-1}\left[  h_{i}^{\left(  \bm  d\right)  }+\sum_{j=2}^{m+n-i+1}%
X_{j}^{\left(  i\right)  }h_{j+i-1}^{\left(  \bm  d\right)  }\right]
+\prod_{j=1}^{k}G_{j};
\end{align*}
For $k=n+m,$%
\begin{align*}
g^{\left(  \bm  d\right)  }\left(  n,m\right)   &  =\left(  -U_{n+m}\right)
^{-1}h_{n+m}^{\left(  \bm  d\right)  }+\sum_{i=1}^{n+m-1}Y_{n+m-i+1}^{\left(
n+m\right)  }\left(  -U_{i}\right)  ^{-1}\\
&  \times\left[  h_{i}^{\left(  \bm  d\right)  }+\sum_{j=2}^{m+n-i+1}%
X_{j}^{\left(  i\right)  }h_{j+i-1}^{\left(  \bm  d\right)  }\right]
+\prod_{j=1}^{n+m}G_{j}.
\end{align*}
\end{The}

\textbf{Proof: }It is seen from (\ref{Eq-2-1}) that we need to compute two
parts: $\mu_{1}\left(  -\mathcal{B}\right)  ^{-1}\bm{e}  _{1}$ and $\left(
-\mathcal{B}\right)  ^{-1}\bm{h}  ^{\left(  \bm      d\right)  }$. For the
first part, we obtain%
\[
\mu_{1}\left(  -\mathcal{B}\right)  ^{-1}\bm{e}  _{1}=\mu_{1}\left(
-U_{1}\right)  ^{-1}\left(  1,Y_{2}^{\left(  2\right)  },Y_{3}^{\left(
3\right)  },\ldots,Y_{n+m}^{\left(  n+m\right)  }\right)  ^{T}.
\]
Since $G_{1}=$ $\left(  -U_{1}\right)  ^{-1}\mu_{1}$ and $Y_{k}^{\left(
k\right)  }=\prod_{j=2}^{k}G_{j}$, we obtain%
\[
\mu_{1}\left(  -\mathcal{B}\right)  ^{-1}\bm{e}  _{1}=\left(  G_{1},G_{2}%
G_{1},G_{3}G_{2}G_{1},\ldots,\prod_{j=1}^{m+n}G_{j}\right)  ^{T}.
\]
For the second part, we have%
\[
\left(  -\mathcal{B}\right)  ^{-1}\bm{h}  ^{\left(  \bm      d\right)
}=\left(
\begin{array}
[c]{l}%
\left(  -U_{1}\right)  ^{-1}\left[  h_{1}^{\left(  \bm  d\right)  }+\sum
_{j=2}^{m+n}X_{j}^{\left(  1\right)  }h_{j}^{\left(  \bm  d\right)  }\right]
\\
\multicolumn{1}{c}{\vdots}\\
\left(  -U_{k}\right)  ^{-1}\left[  h_{k}^{\left(  \bm  d\right)  }+\sum
_{j=2}^{m+n-k+1}X_{j}^{\left(  k\right)  }h_{j+k-1}^{\left(  \bm  d\right)
}\right] \\
+\sum_{i=1}^{k-1}Y_{k-i+1}^{\left(  k\right)  }\left(  -U_{i}\right)
^{-1}\left[  h_{i}^{\left(  \bm  d\right)  }+\sum_{j=2}^{m+n-i+1}%
X_{j}^{\left(  i\right)  }h_{j+i-1}^{\left(  \bm  d\right)  }\right] \\
\multicolumn{1}{c}{\vdots}\\
\left(  -U_{n+m}\right)  ^{-1}h_{n+m}^{\left(  \bm  d\right)  }+\sum
_{i=1}^{n+m-1}Y_{n+m-i+1}^{\left(  n+m\right)  }\left(  -U_{i}\right)  ^{-1}\\
\times\left[  h_{i}^{\left(  \bm  d\right)  }+\sum_{j=2}^{m+n-i+1}%
X_{j}^{\left(  i\right)  }h_{j+i-1}^{\left(  \bm  d\right)  }\right]
\end{array}
\right)  .
\]
Therefore, a simple computation for the vector $\bm      \varphi^{(\bm
d)}=\left(  -\mathcal{B}\right)  ^{-1}\bm{h}  ^{\left(  \bm      d\right)
}+\mu_{1}\left(  -\mathcal{B}\right)  ^{-1}\bm{e}  _{1}$ can derive our
desired results. This completes the proof. \textbf{{\rule{0.08in}{0.08in}}}

\begin{Rem}
(1) Theorem \ref{The:Potential} provides an effective method of solving the
continuous-time Poisson equation $-\mathbf{B}^{\left(  \bm      d\right)
}\bm  g^{(\bm      d)}=\bm      f^{\left(  \bm      d\right)  }-\eta^{\bm
d}\bm  e,$ through an equation transformation%
\begin{equation}
\left\{
\begin{array}
[c]{l}%
g^{\left(  \bm  d\right)  }\left(  0,0\right)  =\Im,\\
\bm  \varphi^{(\bm  d)}=\left(  -\mathcal{B}\right)  ^{-1}\bm{h}  ^{\left(
\bm  d\right)  }+\mu_{1}\left(  -\mathcal{B}\right)  ^{-1}\bm{e}  _{1}\cdot
\Im,
\end{array}
\right.  \label{PoissonS}%
\end{equation}
where $\Im$ is any given constant, and $\left(  -\mathcal{B}\right)  ^{-1}$
can be effectively computed by means of the RG-factorization given in Li
\cite{Li:2010}.

It is necessary to set up a general solution to the system of linear equations
(\ref{PoissonS}). Let $\bm \varphi_{0}^{(\bm d)}$ be the unique solution to
the system of linear equations $-\mathcal{B}\bm \varphi_{0}^{(\bm d)}=\bm{h}
^{\left(  \bm d\right)  }$. Then, the general solution is given by $\bm
\varphi^{(\bm d)}=\bm \varphi_{0}^{(\bm d)}+\mu_{1}\left(  -\mathcal{B}%
\right)  ^{-1}\bm{e} _{1}\cdot\Im.$

(2) To deal with the Poisson equation, some authors (e.g., see Chapter 2 of
Cao \cite{Cao:2007} and Hunter \cite{Hun:1982}) provided a fundamental matrix
method to give a special solution under a constraint condition $\bm
\pi^{\left(  \bm      d\right)  }\bm      g^{(\bm      d)}=\eta^{\bm      d}$,
by which the Poisson equation is well related to the well-known fundamental
matrix (e.g., $\left(  \bm      I-\bm      P+\bm      e\bm      \pi\right)
\bm  g=\bm      f$ or $\left(  -\bm     Q+\bm      e\bm      \pi\right)
\bm      g=\bm  f$). From the fundamental matrix method, our undetermined
constant $\Im$ can be determined by%
\[
\Im=\frac{\eta^{\bm      d}-\bm      \varpi^{\left(  \bm      d\right)
}\left(  -\mathcal{B}\right)  ^{-1}\bm{h}  ^{\left(  \bm      d\right)  }}%
{\pi^{\left(  \bm      d\right)  }\left(  0,0\right)  +\mu_{1}\bm
\varpi^{\left(  \bm  d\right)  }\left(  -\mathcal{B}\right)  ^{-1}\bm{e}
_{1}},
\]
where $\bm      \varpi^{\left(  \bm      d\right)  }$ is a row vector obtained
through omitting the first element of the stationary probability vector $\bm
\pi^{\left(  \bm      d\right)  }$, that is,
\[
\bm      \varpi^{\left(  \bm      d\right)  }=\left(  \pi^{\left(  \bm
d\right)  }\left(  1,0\right)  ,\ldots,\pi^{\left(  \bm      d\right)
}\left(  n,0\right)  ,\pi^{\left(  \bm      d\right)  }\left(  n,1\right)
,\ldots,\pi^{\left(  \bm    d\right)  }\left(  n,m\right)  \right)  .
\]
It is worth noting that in solving the Poisson equation, our RG-factorization
method is superior to the fundamental matrix method because the
RG-factorization given in Li \cite{Li:2010} can easily deal with the inverse
of a high-dimensional transition matrix; while computing the inverse $\left(
\bm      I-\bm      P+\bm      e\bm      \pi\right)  ^{-1}$ or $\left(  -\bm
Q+\bm   e\bm     \pi\right)  ^{-1}$, however, is very difficult for a matrix
$\bm      P$ or $\bm     Q$ with large size, and it also needs to first
compute the stationary probability vector $\bm      \pi$.
\end{Rem}

\section{Impact of Service Price}

In this section, we define a perturbation realization factor of the
policy-based birth-death process, and analyze how the service price impacts on
the perturbation realization factor. Note that the results given in this
section will be utilized for establishing the optimal policy of the
energy-efficient data center in the later section.

For the performance potential vector $\bm \varphi^{(\bm d)}$ under a
constraint condition $g^{\left(  \bm d\right)  }\left(  0,0\right)  =1$, we
define a perturbation realization factor as%
\[
{G}^{\left(  \bm d\right)  }\left(  n,j\right)  \overset{\text{def}}%
{=}g^{\left(  \bm d\right)  }\left(  n,j-1\right)  -g^{\left(  \bm d\right)
}\left(  n,j\right)  ,\text{ }j=1,2,\ldots,m.
\]
It follows from Theorem \ref{The:Potential} that
\begin{align*}
g^{\left(  \bm d\right)  }\left(  n,j-1\right)  = &  \left(  -U_{n+j-1}%
\right)  ^{-1}\left[  h_{n+j-1}^{\left(  \bm d\right)  }+\sum_{k=2}%
^{m-j+2}X_{k}^{\left(  n+j-1\right)  }h_{n+j+k-2}^{\left(  \bm d\right)
}\right]  \\
&  +\sum_{i=1}^{n+j-2}Y_{n+j-i}^{\left(  n+j-1\right)  }\left(  -U_{i}\right)
^{-1}\left[  h_{i}^{\left(  \bm d\right)  }+\sum_{k=2}^{m+n-i+1}X_{k}^{\left(
i\right)  }h_{k+i-1}^{\left(  \bm d\right)  }\right]  +\prod_{k=1}%
^{n+j-1}G_{k}%
\end{align*}
and%
\begin{align*}
g^{\left(  \bm d\right)  }\left(  n,j\right)  = &  \left(  -U_{n+j}\right)
^{-1}\left[  h_{n+j}^{\left(  \bm d\right)  }+\sum_{k=2}^{m-j+1}X_{k}^{\left(
n+j\right)  }h_{n+j+k-1}^{\left(  \bm d\right)  }\right]  \\
&  +\sum_{i=1}^{n+j-1}Y_{n+j-i+1}^{\left(  n+j\right)  }\left(  -U_{i}\right)
^{-1}\left[  h_{i}^{\left(  \bm d\right)  }+\sum_{k=2}^{m+n-i+1}X_{k}^{\left(
i\right)  }h_{k+i-1}^{\left(  \bm d\right)  }\right]  +\prod_{k=1}^{n+j}G_{k}.
\end{align*}
To express the perturbation realization factor ${G}^{\left(  \bm d\right)
}\left(  n,j\right)  $ by means of the service price $R$, we write
\[
A_{0,0}=0,\text{ \ }B_{0,0}=\left(  nP_{1,W}+mP_{2,S}\right)  C_{1}>0;
\]
For $i=1,2,\ldots,n,$
\[
A_{i,0}=i\mu_{1}>0,\text{ \ }B_{i,0}=\left(  nP_{1,W}+mP_{2,S}\right)
C_{1}+iC_{2}^{\left(  1\right)  }>0;
\]
For $j=1,2,\ldots,m$ and $d_{n,j}=0,1,\ldots,m$,
\[
A_{n,j}^{\left(  \bm d\right)  }=n\mu_{1}+\left(  j\wedge d_{n,j}\right)
\mu_{2}>0,
\]%
\begin{align*}
B_{n,j}^{\left(  \bm d\right)  }= &  \left[  nP_{1,W}+d_{n,j}P_{2,W}+\left(
m-d_{n,j}\right)  P_{2,S}\right]  C_{1}+nC_{2}^{\left(  1\right)  }%
+jC_{2}^{\left(  2\right)  }\\
&  +n\mu_{1}C_{3}+\lambda1_{\left\{  i=n,j=m\right\}  }C_{4}>0.
\end{align*}
Then, for $i=0,1,\ldots,n$ and $j=0$,%
\[
f\left(  i,0\right)  =RA_{i,0}-B_{i,0};
\]
For $j=1,2,\ldots,m$,%
\[
f^{\left(  \bm d\right)  }\left(  n,j\right)  =RA_{n,j}^{\left(  \bm d\right)
}-B_{n,j}^{\left(  \bm d\right)  }.
\]
Thus, we obtain%
\begin{align*}
\eta^{\bm d} &  =\bm \pi^{\left(  \bm d\right)  }\bm f^{\left(  \bm d\right)
}\\
&  =\sum_{i=0}^{n}\pi\left(  i,0\right)  f\left(  i,0\right)  +\sum_{j=1}%
^{m}\pi\left(  n,j\right)  f^{\left(  \bm d\right)  }\left(  n,j\right)  \\
&  =R{D}^{\left(  \bm d\right)  }-{F}^{\left(  \bm d\right)  },
\end{align*}
where%
\[
{D}^{\left(  \bm d\right)  }=\sum_{i=0}^{n}\pi\left(  i,0\right)  A_{i,0}%
+\sum_{j=1}^{m}\pi\left(  n,j\right)  A_{n,j}^{\left(  \bm d\right)  }>0,
\]
and%
\[
{F}^{\left(  \bm d\right)  }=\sum_{i=0}^{n}\pi\left(  i,0\right)  B_{i,0}%
+\sum_{j=1}^{m}\pi\left(  n,j\right)  B_{n,j}^{\left(  \bm d\right)  }>0.
\]
Then,%
\[
\bm{h} ^{\left(  \bm d\right)  }=\left(
\begin{array}
[c]{c}%
h_{1}^{\left(  \bm d\right)  }\\
\vdots\\
h_{n}^{\left(  \bm d\right)  }\\
h_{n+1}^{\left(  \bm d\right)  }\\
\vdots\\
h_{n+m}^{\left(  \bm d\right)  }%
\end{array}
\right)  =\left(
\begin{array}
[c]{c}%
f\left(  1,0\right)  -\eta^{\bm d}\\
\vdots\\
f\left(  n,0\right)  -\eta^{\bm d}\\
f^{\left(  \bm d\right)  }\left(  n,1\right)  -\eta^{\bm d}\\
\vdots\\
f^{\left(  \bm d\right)  }\left(  n,m\right)  -\eta^{\bm d}%
\end{array}
\right)  =\left(
\begin{array}
[c]{c}%
R\left[  A_{1,0}-{D}^{\left(  \bm d\right)  }\right]  -\left[  B_{1,0}%
-{F}^{\left(  \bm d\right)  }\right]  \\
\vdots\\
R\left[  A_{n,0}-{D}^{\left(  \bm d\right)  }\right]  -\left[  B_{n,0}%
-{F}^{\left(  \bm d\right)  }\right]  \\
R\left[  A_{n,1}^{\left(  \bm d\right)  }-{D}^{\left(  \bm d\right)  }\right]
-\left[  B_{n,1}^{\left(  \bm d\right)  }-{F}^{\left(  \bm d\right)  }\right]
\\
\vdots\\
R\left[  A_{n,m}^{\left(  \bm d\right)  }-{D}^{\left(  \bm d\right)  }\right]
-\left[  B_{n,m}^{\left(  \bm d\right)  }-{F}^{\left(  \bm d\right)
}\right]
\end{array}
\right)  .
\]

If a job finishes its service at a server and leaves this system immediately,
then the data center can obtain a fixed revenue $R$ from each job. Obviously,
$R$ is the service price provided by the data center. Now, we study the
influence of the service price $R$ on the perturbation realization factor
${G}^{\left(  \bm d\right)  }\left(  n,j\right)  $. Note that all the numbers
$\left(  -U_{k}\right)  ^{-1}$, $X_{j}^{\left(  k\right)  }$, $Y_{j}^{\left(
k\right)  }$ and $G_{j}$ are positive and are independent of the service price
$R$, while all the numbers $h_{j}^{\left(  \bm d\right)  }$ are the linear
functions of $R$.

We write%
\begin{align*}
{W}_{n,j}^{\left(  \bm d\right)  }= &  \left(  -U_{n+j}\right)  ^{-1}\left\{
\left[  A_{n,j}^{\left(  \bm d\right)  }-{D}^{\left(  \bm d\right)  }\right]
+\sum_{k=2}^{m-j+1}X_{k}^{\left(  n+j\right)  }\left[  A_{n,j+k-1}^{\left(
\bm d\right)  }-{D}^{\left(  \bm d\right)  }\right]  \right\}  \\
&  +\sum_{i=1}^{n}Y_{n+j-i+1}^{\left(  n+j\right)  }\left(  -U_{i}\right)
^{-1}\left[  A_{i,0}-{D}^{\left(  \bm d\right)  }\right]  +\sum_{i=n+1}%
^{n+j-1}Y_{n+j-i+1}^{\left(  n+j\right)  }\left(  -U_{i}\right)  ^{-1}\left[
A_{n,i-n}^{\left(  \bm d\right)  }-{D}^{\left(  \bm d\right)  }\right]  \\
&  +\sum_{3\leq k+i\leq n+1}\sum_{i=1}^{n+j-1}Y_{n+j-i+1}^{\left(  n+j\right)
}\left(  -U_{i}\right)  ^{-1}\sum_{k=2}^{m+n-i+1}X_{k}^{\left(  i\right)
}\left[  A_{k+i-1,0}-{D}^{\left(  \bm d\right)  }\right]  \\
&  +\sum_{n+2\leq k+i\leq n+m+1}\sum_{i=1}^{n+j-1}Y_{n+j-i+1}^{\left(
n+j\right)  }\left(  -U_{i}\right)  ^{-1}\sum_{k=2}^{m+n-i+1}X_{k}^{\left(
i\right)  }\left[  A_{n,k+i-1-n}^{\left(  \bm d\right)  }-{D}^{\left(  \bm
d\right)  }\right]
\end{align*}
and%
\begin{align*}
{V}_{n,j}^{\left(  \bm d\right)  }= &  \left(  -U_{n+j}\right)  ^{-1}\left\{
\left[  B_{n,j}^{\left(  \bm d\right)  }-{F}^{\left(  \bm d\right)  }\right]
+\sum_{k=2}^{m-j+1}X_{k}^{\left(  n+j\right)  }\left[  B_{n,j+k-1}^{\left(
\bm d\right)  }-{F}^{\left(  \bm d\right)  }\right]  \right\}  \\
&  +\sum_{i=1}^{n}Y_{n+j-i+1}^{\left(  n+j\right)  }\left(  -U_{i}\right)
^{-1}\left[  B_{i,0}-{F}^{\left(  \bm d\right)  }\right]  +\sum_{i=n+1}%
^{n+j-1}Y_{n+j-i+1}^{\left(  n+j\right)  }\left(  -U_{i}\right)  ^{-1}\left[
B_{n,i-n}^{\left(  \bm d\right)  }-{F}^{\left(  \bm d\right)  }\right]  \\
&  +\sum_{3\leq k+i\leq n+1}\sum_{i=1}^{n+j-1}Y_{n+j-i+1}^{\left(  n+j\right)
}\left(  -U_{i}\right)  ^{-1}\sum_{k=2}^{m+n-i+1}X_{k}^{\left(  i\right)
}\left[  B_{k+i-1,0}-{F}^{\left(  \bm d\right)  }\right]  \\
&  +\sum_{n+2\leq k+i\leq n+m+1}\sum_{i=1}^{n+j-1}Y_{n+j-i+1}^{\left(
n+j\right)  }\left(  -U_{i}\right)  ^{-1}\sum_{k=2}^{m+n-i+1}X_{k}^{\left(
i\right)  }\left[  B_{n,k+i-1-n}^{\left(  \bm d\right)  }-{F}^{\left(  \bm
d\right)  }\right]  ,
\end{align*}
then, we obtain that for $j=1,2,\ldots,m$,%
\begin{align}
{G}^{\left(  \bm d\right)  }\left(  n,j\right)   &  =g^{\left(  \bm d\right)
}\left(  n,j-1\right)  -g^{\left(  \bm d\right)  }\left(  n,j\right)
\nonumber\\
&  =R\left[  {W}_{n,j-1}^{\left(  \bm d\right)  }-{W}_{n,j}^{\left(  \bm
d\right)  }\right]  -\left[  {V}_{n,j-1}^{\left(  \bm d\right)  }-{V}%
_{n,j}^{\left(  \bm d\right)  }\right]  +\left(  1-G_{n+j}\right)  \prod
_{k=1}^{n+j-1}G_{k}.\label{Eq-9}%
\end{align}

We can see that ${G}^{\left(  \bm      d\right)  }\left(  n,j\right)  $
quantifies the difference among two adjacent performance potentials
$g^{\left(  \bm     d\right)  }\left(  n,j\right)  $ and $g^{\left(  \bm
d\right)  }\left(  n,j-1\right)  $. It measures the long-run effect on the
average profit of the data center when the system state is changed from
$(n,j-1)$ to $(n,j)$, which indicates the occurrence of a service completion
event. From later discussion in Section 6, we will see that ${G}^{\left(
\bm      d\right)  }\left(  n,j\right)  $\ plays a fundamental role in the
performance optimization of data centers and the sign of ${G}^{\left(  \bm
d\right)  }\left(  n,j\right)  +c$ directly determines the selection of
decision actions, as shown in (\ref{eq-19})\ later, where $c$ is defined as%
\begin{equation}
c=R-\frac{\left(  P_{2,W}-P_{2,S}\right)  C_{1}}{\mu_{2}}. \label{Eq-10}%
\end{equation}
To this end, we analyze how the service price impacts on ${G}^{\left(  \bm
d\right)  }\left(  n,j\right)  +c$ as follows.

Substituting (\ref{Eq-9}) into the linear equation ${G}^{\left(  \bm
d\right)  }\left(  n,j\right)  +c=0$, we obtain
\begin{equation}
R\left[  {W}_{n,j-1}^{\left(  \bm  d\right)  }-{W}_{n,j}^{\left(  \bm
d\right)  }\right]  -\left[  {V}_{n,j-1}^{\left(  \bm  d\right)  }-{V}%
_{n,j}^{\left(  \bm  d\right)  }\right]  \mathbf{+}\left(  1-G_{n+j}\right)
\prod_{k=1}^{n+j-1}G_{k}+c=0. \label{Eq-4}%
\end{equation}
Substituting (\ref{Eq-10}) into the above equation, we obtain that the unique
solution of the price $R$ in (\ref{Eq-4}) is given by%
\begin{equation}
\Re_{n,j}^{\left(  \bm  d\right)  }=\frac{\left[  {V}_{n,j-1}^{\left(  \bm
d\right)  }-{V}_{n,j}^{\left(  \bm  d\right)  }\right]  \mathbf{-}\left(
1-G_{n+j}\right)  \prod_{k=1}^{n+j-1}G_{k}+\frac{\left(  P_{2,W}%
-P_{2,S}\right)  C_{1}}{\mu_{2}}}{1+\left[  {W}_{n,j-1}^{\left(  \bm
d\right)  }-{W}_{n,j}^{\left(  \bm  d\right)  }\right]  }. \label{Eq-5}%
\end{equation}
It is easy to see that (a) if $R\geq\Re_{n,j}^{\left(  \bm      d\right)  }$,
then ${G}^{\left(  \bm      d\right)  }\left(  n,j\right)  +{c}\geq0$; and (b)
if $R\leq\Re_{n,j}^{\left(  \bm      d\right)  }$, then ${G}^{\left(  \bm
d\right)  }\left(  n,j\right)  +{c}\leq0$.

In the energy-efficient data center, we define two critical values, related to
the service price, as%
\begin{equation}
{R}_{H}=\max_{\bm  d\in\mathcal{D}}\left\{  0,\Re_{n,1}^{\left(  \bm
d\right)  },\Re_{n,2}^{\left(  \bm  d\right)  },\ldots,\Re_{n,m}^{\left(  \bm
d\right)  } \right\}  \label{Eq-6}%
\end{equation}
and%
\begin{equation}
{R}_{L}=\min_{\bm  d\in\mathcal{D}}\left\{  \Re_{n,1}^{\left(  \bm  d\right)
},\Re_{n,2}^{\left(  \bm  d\right)  },\ldots,\Re_{n,m}^{\left(  \bm  d\right)
} \right\}  . \label{Eq-6-1}%
\end{equation}

The following proposition uses the two critical values related to the service
price to provide a key condition whose purpose is to establish a
sensitivity-based optimization framework of the energy-efficient data center
in our later study. Also, this proposition will be useful in the next section
for studying the monotonicity of the energy-efficient policies.

\begin{Pro}
\label{Pro:servicep}(1) If $R\geq{R}_{H}$, then for any $\bm  d\in\mathcal{D}$
and for each $j=1,2,\ldots,m$, we have%
\begin{equation}
{G}^{\left(  \bm  d\right)  }\left(  n,j\right)  +{c}\geq0\text{.}
\label{Eq-7}%
\end{equation}

(2) If $0\leq R\leq{R}_{L}$, then for any $\bm                  d\in
\mathcal{D}$ and for each $j=1,2,\ldots,m$, we have%
\begin{equation}
{G}^{\left(  \bm  d\right)  }\left(  n,j\right)  +{c}\leq0\text{.}
\label{Eq-7-1}%
\end{equation}
\end{Pro}

\textbf{Proof: }(1) For any $\bm            d\in\mathcal{D}$ and for each
$j=1,2,\ldots,m$, since $R\geq{R}_{H}$ and ${R}_{H}=\max_{\bm   d\in
\mathcal{D}}\left\{  0,\Re_{n,1}^{\left(  \bm                  d\right)  }%
,\Re_{n,2}^{\left(  \bm           d\right)  },\ldots,\Re_{n,m}^{\left(  \bm
d\right)  }\right\}  $, we have%
\[
R\geq\Re_{n,j}^{\left(  \bm                  d\right)  },
\]
which clearly makes that ${G}^{\left(  \bm                  d\right)  }\left(
n,j\right)  +{c}\geq0$.

(2) For any $\bm         d\in\mathcal{D}$ and for each $j=1,2,\ldots,m$, if
$0\leq R\leq{R}_{L}$, we have%
\[
R\leq\Re_{n,j}^{\left(  \bm         d\right)  },
\]
this gives that ${G}^{\left(  \bm         d\right)  }\left(  n,j\right)
+{c}\leq0$. This completes the proof. \textbf{{\rule{0.08in}{0.08in}}}

\section{Monotonicity and Optimality}

In this section, we use the Poisson equation to derive a useful performance
difference equation, and discuss the monotonicity and optimality of the
long-run average profit of the energy-efficient data center with respect to
the policies. Based on this, we give the optimal energy-efficient policy under
some restrained service prices.

For any given policy $\bm                                     d\in\mathcal{D}%
$, the policy-based continuous-time birth-death process $\{\mathbf{X}^{\left(
\bm                      d\right)  }(t):t\geq0\}$ with infinitesimal generator
$\mathbf{B}^{\left(  \bm                   d\right)  }$ given in
(\ref{eq-60})\ is irreducible, aperiodic and positive recurrent, hence the
long-run average profit of the data center is given by
\[
\eta^{\bm                                     d}=\bm            \pi^{\left(
\bm   d\right)  }\bm                            f^{\left(  \bm  d\right)  },
\]
and the Poisson equation is written as%
\[
\mathbf{B}^{\left(  \bm                                     d\right)  }\bm
g^{(\bm        d)} = \eta^{\bm                                d}
\bm                                    e - \bm           f^{\left(  \bm
d\right)  }.
\]
With a similar role played by state $\left(  i,j\right)  $, it is seen from
(\ref{eq-60}) that the policy $\bm                                     d$
directly affects not only the elements of the infinitesimal generator
$\mathbf{B}^{\left(  \bm             d\right)  }$ but also the reward function
$\bm     f^{\left(  \bm              d\right)  }$. That is, if the policy
$\bm  d$ changes, then the infinitesimal generator $\mathbf{B}^{\left(  \bm
d\right)  }$ and the reward function $\bm    f^{\left(  \bm     d\right)  }$
will have their corresponding changes. To express such a change
mathematically, we take two different policies $\bm       d$ and $\bm
d^{\prime}$, both of which correspond to their infinitesimal generators
$\mathbf{B}^{\left(  \bm   d\right)  }$ and $\mathbf{B}^{\left(  \bm
d^{\prime}\right)  }$, and to their reward functions $\bm       f^{\left(
\bm    d\right)  }$ and $\bm   f^{\left(  \bm                d^{\prime
}\right)  }$.

The following lemma provides a useful equation for the difference $\eta^{\bm
d^{\prime}}-\eta^{\bm     d}$ of the long-run average performances $\eta^{\bm
d}$ and $\eta^{\bm     d^{\prime}}$ for any two policies $\bm     d,\bm
d^{\prime}\in\mathcal{D}$. The performance difference equation plays a key
role in the sensitivity-based optimization theory. Note that the performance
difference equation was given in Cao \cite{Cao:2007}, while here we restate it
with some simple discussion, for convenience of readers.

\begin{Lem}
\label{Lem:diff}For any two policies $\bm      d,\bm       d^{\prime}%
\in\mathcal{D}$, we have
\begin{equation}
\eta^{\bm  d^{\prime}}-\eta^{\bm  d}=\bm  \pi^{\left(  \bm  d^{\prime}\right)
}\left[  \left(  \mathbf{B}^{\left(  \bm  d^{\prime}\right)  }-\mathbf{B}%
^{\left(  \bm  d\right)  }\right)  \bm  g^{(\bm  d)}\mathbf{+}\left(  \bm
f^{\left(  \bm  d^{\prime}\right)  }-\bm  f^{\left(  \bm  d\right)  }\right)
\right]  . \label{eq-16}%
\end{equation}
\end{Lem}

\textbf{Proof: }Note that $\bm     \pi^{\left(  \bm     d^{\prime}\right)
}\mathbf{B}^{\left(  \bm     d^{\prime}\right)  }=\mathbf{0}$, $\mathbf{B}%
^{\left(  \bm     d\right)  }\bm     g^{(\bm     d)}=\eta^{\bm     d}\bm
e-\bm  f^{\left(  \bm     d\right)  }$, $\bm     \pi^{\left(  \bm
d^{\prime}\right)  }\bm    e=1$, we compute%
\begin{align*}
&  \bm  \pi^{\left(  \bm  d^{\prime}\right)  }\left[  \left(  \mathbf{B}%
^{\left(  \bm  d^{\prime}\right)  }-\mathbf{B}^{\left(  \bm  d\right)
}\right)  \bm  g^{(\bm  d)}\mathbf{+}\left(  \bm  f^{\left(  \bm  d^{\prime
}\right)  }-\bm  f^{\left(  \bm  d\right)  }\right)  \right] \\
&  =-\bm  \pi^{\left(  \bm  d^{\prime}\right)  }\cdot\mathbf{B}^{\left(  \bm
d\right)  }\bm  g^{(\bm  d)}+\bm  \pi^{\left(  \bm  d^{\prime}\right)  }\bm
f^{\left(  \bm  d^{\prime}\right)  }-\bm  \pi^{\left(  \bm  d^{\prime}\right)
}\bm  f^{\left(  \bm  d\right)  }\\
&  =-\bm  \pi^{\left(  \bm  d^{\prime}\right)  }\left[  \eta^{\bm  d}\bm
e-\bm  f^{\left(  \bm  d\right)  }\right]  +\eta^{\bm  d^{\prime}}-\bm
\pi^{\left(  \bm  d^{\prime}\right)  }\bm  f^{\left(  \bm  d\right)  }\\
&  =\eta^{\bm  d^{\prime}}-\eta^{\bm  d}.
\end{align*}
This completes the proof. \textbf{{\rule{0.08in}{0.08in}}}

Now, we describe the first role played by the performance difference, in which
we set up a partial order relation in the policy set $\mathcal{D}$ so that the
optimal policy in $\mathcal{D}$ can be found numerically. Based on the
performance difference $\eta^{\bm  d^{\prime}}-\eta^{\bm               d}$ for
any two policies $\bm  d,\bm   d^{\prime}\in\mathcal{D}$, we can set up a
partial order in the policy set $\mathcal{D}$ as follows. We write
$\bm            d^{\prime}\succ$ $\bm  d$ if $\eta
^{\bm                         d^{\prime}}>\eta^{\bm              d}$; $\bm
d^{\prime}\approx\bm            d$ if $\eta^{\bm              d^{\prime}}%
=\eta^{\bm    d}$; $\bm                                   d^{\prime}%
\prec\bm                d$ if $\eta^{\bm        d^{\prime}}<\eta
^{\bm                                     d}$. Also, we write $\bm
d^{\prime}\succeq\bm    d$ if $\eta^{\bm         d^{\prime}}\geq
\eta^{\bm                         d}$; $\bm    d^{\prime}\preceq
\bm                                     d$ if $\eta^{\bm    d^{\prime}}%
\leq\eta^{\bm                   d}$. By using this partial order, our research
target is to find an optimal policy $\bm               d^{\ast}\in\mathcal{D}$
such that $\bm        d^{\ast}\succeq\bm
d$ for any policy $\bm        d\in\mathcal{D}$, or%
\[
\bm                                     d^{\ast}=\underset{\bm
d\in\mathcal{D}}{\arg\max}\left\{  \eta^{\bm        d}\right\}  .
\]

Note that the policy set $\mathcal{D}$ and the state space $\bm     \Omega$
are both finite, thus an enumeration method is feasible for finding the
optimal energy-efficient policy $\bm     d^{\ast}$ in the policy set
$\mathcal{D}$. Since $\bm     d=\left(  0,0,\ldots,0;d_{n,1},d_{n,2}%
,\ldots,d_{n,m}\right)  $ and $d_{n,j}\in\left\{  0,1,\ldots,m\right\}  $, it
is seen that the policy set $\mathcal{D}$ contains $\left(  m+1\right)  ^{m}$
elements so that the enumeration method will require a huge computation
workload. However, our following work can greatly reduce the optimization
complexity by means of the sensitivity-based optimization theory.

Now, we discuss the monotonicity of the long-run average profit $\eta^{\bm d}$
with respect to a decision element $d_{n,j}$ of any policy $\bm d\in
\mathcal{D}$, for $d_{n,j}=0,1,\ldots,m$. This result is derived by the
following two theorems, in which we show that for any policy $\bm
d\in\mathcal{D}$ and for each $j=1,2,\ldots,m$, the long-run average profit
$\eta^{\bm d}$ is unimodal with respect to each decision element $d_{n,j}%
\in\left\{  0,1,\ldots,m\right\}  $.

\begin{The}
\label{The:right}For any policy $\bm     d\in\mathcal{D}$ and for each
$j=1,2,\ldots,m$, the long-run average profit $\eta^{\bm     d}$ is linearly
decreasing with respect to each decision element $d_{n,j}$, where $d_{n,j}%
\in\left\{  j,j+1,\ldots,m\right\}  $.
\end{The}

\textbf{Proof: }For each $j=1,2,\ldots,m$, we consider two interrelated
policies $\bm d,\bm d^{\prime}\in\mathcal{D}$ as follows.%
\begin{align*}
\bm d &  =\left(  0,0,\ldots,0;d_{n,1},d_{n,2},\ldots,d_{n,j-1}\underline
{,d_{n,j},}d_{n,j+1},\ldots,d_{n,m}\right)  ,\\
\bm d^{\prime} &  =\left(  0,0,\ldots,0;d_{n,1},d_{n,2},\ldots,d_{n,j-1}%
\underline{,j,}d_{n,j+1},\ldots,d_{n,m}\right)  ,
\end{align*}
where $d_{n,j}>j$. It is seen that the two policies $\bm d,\bm d^{\prime}$
have one difference only between their corresponding decision elements
$d_{n,j}$ and $j$. In this case, it is seen from Theorem 1 that $\mathbf{B}%
^{\left(  \bm d\right)  }=\mathbf{B}^{\left(  \bm d^{\prime}\right)  }$ and
$\bm \pi^{\left(  \bm d\right)  }=\bm \pi^{\left(  \bm d^{\prime}\right)  }$.
Also, it is easy to check from (\ref{eq-51}) to (\ref{eq-33}) that%
\[
\bm f^{\left(  \bm d\right)  }-\bm f^{\left(  \bm d^{\prime}\right)  }=\left(
0,0,\ldots,0;0,0,\ldots,0\underline{,-\left(  d_{n,j}-j\right)  \left(
P_{2,W}-P_{2,S}\right)  C_{1},}0,\ldots,0\right)  ^{T}.
\]
Thus, it follows from Lemma \ref{Lem:diff} that%
\begin{align*}
\eta^{\bm d}-\eta^{\bm d^{\prime}} &  =\bm \pi^{\left(  \bm d\right)  }\left[
\left(  \mathbf{B}^{\left(  \bm d\right)  }-\mathbf{B}^{\left(  \bm d^{\prime
}\right)  }\right)  \bm g^{(\bm d^{\prime})}\mathbf{+}\left(  \bm f^{\left(
\bm d\right)  }-\bm f^{\left(  \bm d^{\prime}\right)  }\right)  \right]  \\
&  =-\pi^{\left(  \bm d\right)  }\left(  n,j\right)  \left(  d_{n,j}-j\right)
\left(  P_{2,W}-P_{2,S}\right)  C_{1}%
\end{align*}
or%
\begin{equation}
\eta^{\bm d}=\eta^{\bm d^{\prime}}-\pi^{\left(  \bm d\right)  }\left(
n,j\right)  \left(  d_{n,j}-j\right)  \left(  P_{2,W}-P_{2,S}\right)
C_{1}.\label{Eq-8}%
\end{equation}
Since $\bm \pi^{\left(  \bm d\right)  }=\bm \pi^{\left(  \bm d^{\prime
}\right)  }$ by Theorem 1, it is easy to see that $\pi^{\left(  \bm d\right)
}\left(  n,j\right)  =\pi^{\left(  \bm d^{\prime}\right)  }\left(  n,j\right)
$\ can be determined by $d_{n,j}^{\prime}=j.$ This indicates that
$\pi^{\left(  \bm d\right)  }\left(  n,j\right)  $ is irrelevant to the
decision element $d_{n,j}.$ Again, note that $\eta^{\bm d^{\prime}}$ is
irrelevant to the decision element $d_{n,j}$, and $P_{2,W}-P_{2,S}$ and
$C_{1}$ are two positive constants. Thus, it is easy to see from (\ref{Eq-8})
that the long-run average profit $\eta^{\bm d}$ is linearly decreasing with
respect to each decision element $d_{n,j}\in\left\{  j,j+1,\ldots,m\right\}
$. This completes the proof. \textbf{{\rule{0.08in}{0.08in}}}

In what follows we discuss the left half part of the unimodal structure of the
long-run average profit $\eta^{\bm d}$ with respect to each decision element
$d_{n,j}\in\left\{  0,1,\ldots,j\right\}  $. Compared to the analysis of its
right half part in Theorem 3, our discussion for the left half part is a
little bit complicated.

Let the optimal energy-efficient policy $\bm     d^{\ast}=\underset{\bm
d\in\mathcal{D}}{\arg\max}\left\{  \eta^{\bm     d}\right\}  $ be%
\[
\bm     d^{\ast}=\left(  0,0,\ldots,0;d_{n,1}^{\ast},d_{n,2}^{\ast}%
,\ldots,d_{n,m}^{\ast}\right)  .
\]
Then, it is seen from Theorem \ref{The:right} that
\begin{gather*}
d_{n,1}^{\ast}\in\left\{  0,1\right\}  ;\\
\vdots\\
d_{n,j}^{\ast}\in\left\{  0,1,\ldots,j\right\}  ;\\
\vdots\\
d_{n,m}^{\ast}\in\left\{  0,1,\ldots,m\right\}  .
\end{gather*}
Thus, Theorem \ref{The:right} makes the area of finding the optimal
energy-efficient policy $\bm     d^{\ast}$ from a large set $\left\{
0,1,\ldots,m\right\}  ^{m}$ to a shrunken area $\left\{  0,1\right\}
\times\left\{  0,1,2\right\}  \times\cdots\times\left\{  0,1,\ldots,m\right\}
$.

To find the optimal energy-efficient policy $\bm d^{\ast}$, we consider two
energy-efficient policies with an interrelated structure as follows.%
\begin{align*}
\bm d &  =\left(  0,0,\ldots,0;d_{n,1},\ldots,d_{n,j-1}\underline{,d_{n,j}%
,}d_{n,j+1},\ldots,d_{n,m}\right)  ,\\
\bm d^{\prime} &  =\left(  0,0,\ldots,0;d_{n,1},\ldots,d_{n,j-1}%
\underline{,d_{n,j}^{\prime},}d_{n,j+1},\ldots,d_{n,m}\right)  ,
\end{align*}
where $d_{n,j}^{\prime}>d_{n,j}$, and $d_{n,j},d_{n,j}^{\prime}\in\left\{
1,2,\ldots,j\right\}  $. It is easy to check from (\ref{eq-60}) that
\begin{equation}
\mathbf{B}^{\left(  \bm d^{\prime}\right)  }-\mathbf{B}^{\left(  \bm d\right)
}\mathbf{=}\left(
\begin{array}
[c]{ccccccc}%
0 &  &  &  &  &  & \\
0 & \ddots &  &  &  &  & \\
& \ddots & 0 &  &  &  & \\
&  & \left(  d_{n,j}^{\prime}-d_{n,j}\right)  \mu_{2} & -\left(
d_{n,j}^{\prime}-d_{n,j}\right)  \mu_{2} &  &  & \\
&  &  & 0 & 0 &  & \\
&  &  &  & \ddots & \ddots & \\
&  &  &  &  & 0 & 0
\end{array}
\right)  .\label{eq-17}%
\end{equation}
On the other hand, from the reward functions given in (\ref{eq-51})$\ $to
(\ref{eq-53}), it is seen that for $j=1,2,\ldots,m$, and $d_{n,j}%
=0,1,\ldots,m$,%
\begin{align*}
f^{\left(  \bm d\right)  }\left(  n,j\right)  = &  \left[  R\mu_{2}-\left(
P_{2,W}-P_{2,S}\right)  C_{1}\right]  d_{n,j}\\
&  +Rn\mu_{1}-\left(  nP_{1,W}-mP_{2,S}\right)  C_{1}-\left[  nC_{2}^{\left(
1\right)  }+jC_{2}^{\left(  2\right)  }\right]  -n\mu_{1}1_{\{j>0\}}%
C_{3}-\lambda1_{\left\{  i=n,j=m\right\}  }C_{4}%
\end{align*}
and%
\begin{align*}
f^{\left(  \bm d^{\prime}\right)  }\left(  n,j\right)  = &  \left[  R\mu
_{2}-\left(  P_{2,W}-P_{2,S}\right)  C_{1}\right]  d_{n,j}^{\prime}\\
&  +Rn\mu_{1}-\left(  nP_{1,W}-mP_{2,S}\right)  C_{1}-\left[  nC_{2}^{\left(
1\right)  }+jC_{2}^{\left(  2\right)  }\right]  -n\mu_{1}1_{\{j>0\}}%
C_{3}-\lambda1_{\left\{  i=n,j=m\right\}  }C_{4}.
\end{align*}
Hence, we have
\begin{equation}
\bm f^{\left(  \bm d^{\prime}\right)  }-\bm f^{\left(  \bm d\right)  }=\left(
0,0,\ldots,0,\mu_{2}{c}\left(  d_{n,j}^{\prime}-d_{n,j}\right)  ,0,\ldots
,0\right)  ^{T}.\label{eq-18}%
\end{equation}
We write
\[
\eta^{\bm d}|_{d_{n,j}=k}=\bm \pi^{\left(  \bm d\right)  }|_{d_{n,j}=k}%
\cdot\bm f^{\left(  \bm d\right)  }|_{d_{n,j}=k}.
\]
The following theorem discusses the left half part of the unimodal structure
of the long-run average profit $\eta^{\bm d}$ with respect to each decision
element $d_{n,j}\in\left\{  0,1,\ldots,m\right\}  $.

\begin{The}
\label{The:left}If $R\geq{R}_{H}$, then for any policy $\bm     d\in
\mathcal{D}$ and for each $j=1,2,\ldots,m$, the long-run average profit
$\eta^{\bm     d}$ is strictly monotone increasing with respect to each
decision element $d_{n,j}$, where $d_{n,j}\in\left\{  0,1,\ldots,j\right\}  $.
\end{The}

\textbf{Proof: }For each $j=1,2,\ldots,m$, we consider two energy-efficient
policies with an interrelated structure as follows.%
\begin{align*}
\bm  d  &  =\left(  0,0,\ldots,0;d_{n,1},\ldots,d_{n,j-1}\underline{,d_{n,j}%
,}d_{n,j+1},\ldots,d_{n,m}\right)  ,\\
\bm  d^{\prime}  &  =\left(  0,0,\ldots,0;d_{n,1},\ldots,d_{n,j-1}%
\underline{,d_{n,j}^{\prime},}d_{n,j+1},\ldots,d_{n,m}\right)  ,
\end{align*}
where $d_{n,j}^{\prime}>d_{n,j}$, and $d_{n,j},d_{n,j}^{\prime}\in\left\{
0,1,\ldots,j\right\}  $. Applying Lemma \ref{Lem:diff}, it follows from
(\ref{eq-17})$\ $and (\ref{eq-18}) that%
\begin{align}
\eta^{\bm  d^{\prime}}-\eta^{\bm  d}  &  =\bm  \pi^{\left(  \bm  d^{\prime
}\right)  }\left[  \left(  \mathbf{B}^{\left(  \bm  d^{\prime}\right)
}-\mathbf{B}^{\left(  \bm  d\right)  }\right)  \bm  g^{(\bm  d)}%
\mathbf{+}\left(  \bm  f^{\left(  \bm  d^{\prime}\right)  }-\bm  f^{\left(
\bm  d\right)  }\right)  \right] \nonumber\\
&  =\mu_{2}\pi^{\left(  \bm  d^{\prime}\right)  }\left(  n,j\right)  \left(
d_{n,j}^{\prime}-d_{n,j}\right)  \left[  g^{\left(  \bm  d\right)  }\left(
n,j-1\right)  -g^{\left(  \bm  d\right)  }\left(  n,j\right)  +{c}\right]
\nonumber\\
&  =\mu_{2}\pi^{\left(  \bm  d^{\prime}\right)  }\left(  n,j\right)  \left(
d_{n,j}^{\prime}-d_{n,j}\right)  \left[  {G}^{\left(  \bm  d\right)  }\left(
n,j\right)  +{c}\right]  , \label{eq-19}%
\end{align}
where ${G}^{\left(  \bm     d\right)  }\left(  n,j\right)  =g^{\left(  \bm
d\right)  }\left(  n,j-1\right)  -g^{\left(  \bm     d\right)  }\left(
n,j\right)  $. If $R\geq{R}_{H}$, then it is seen from Proposition
\ref{Pro:servicep} that ${G}^{\left(  \bm     d\right)  }\left(  n,j\right)
+{c\geq0}$. Thus, we obtain that for the two policies $\bm     d,\bm
d^{\prime}\in\mathcal{D}$ with $d_{n,j}^{\prime}>d_{n,j}$ and $d_{n,j}%
,d_{n,j}^{\prime}\in\left\{  0,1,\ldots,j\right\}  $,%
\[
\eta^{\bm     d^{\prime}}>\eta^{\bm     d}.
\]
This shows that%
\[
\eta^{\bm     d}|_{d_{n,j}=1}<\eta^{\bm     d}|_{d_{n,j}=2}<\cdots<\eta^{\bm
d}|_{d_{n,j}=m-1}<\eta^{\bm     d}|_{d_{n,j}=m}.
\]
This completes the proof. \textbf{{\rule{0.08in}{0.08in}}}

When $R\geq{R}_{H}$, we use Figure 2 to provide an intuitive summary for the
main results given in Theorems \ref{The:right} and \ref{The:left}. In the
right half part of Figure 2,%
\[
\eta^{\bm    d}=\eta^{\bm    d^{\prime}}-\pi^{\left(  \bm    d\right)
}\left(  n,j\right)  \left(  d_{n,j}-j\right)  \left(  P_{2,W}-P_{2,S}\right)
C_{1}%
\]
shows that $\eta^{\bm    d}$ is a linear function of the decision element
$d_{n,j}$. By contrast, in the right half part of Figure 2, we need to first
introduce a restrictive condition: $R\geq{R}_{H}$, under which%
\[
\eta^{\bm    d^{\prime}}-\eta^{\bm    d}=\mu_{2}\pi^{\left(  \bm    d^{\prime
}\right)  }\left(  n,j\right)  \left(  d_{n,j}^{\prime}-d_{n,j}\right)
\left[  {G}^{\left(  \bm    d\right)  }\left(  n,j\right)  +{c}\right]  .
\]
Let $d_{n,j}^{\prime}=j$. Then,%
\[
\eta^{\bm    d}=\eta^{\bm    d^{\prime}}-\mu_{2}\pi^{\left(  \bm    d^{\prime
}\right)  }\left(  n,j\right)  \left(  j-d_{n,j}\right)  \left[  {G}^{\left(
\bm    d\right)  }\left(  n,j\right)  +{c}\right]  .
\]
Since ${G}^{\left(  \bm    d\right)  }\left(  n,j\right)  $ depends on the
decision element $d_{n,j}$, it is clear that $\eta^{\bm    d}$ is a nonlinear
function of the decision element $d_{n,j}$.

\begin{figure}[th]
\centering          \includegraphics[width=12cm]{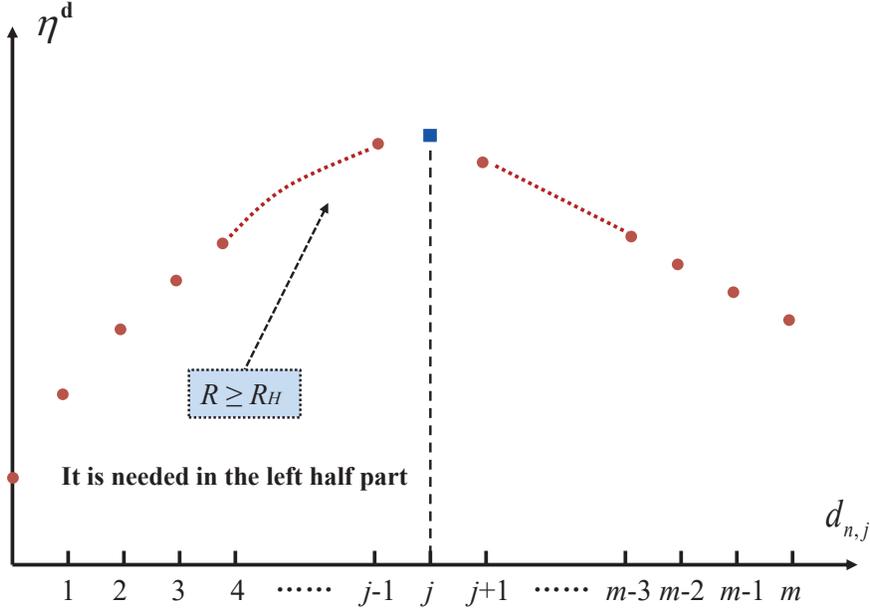}  \caption{The
unimodal structure of the long-run average profit.}%
\label{figure:fig-2}%
\end{figure}

\begin{The}
\label{Cor:left}If $0\leq R\leq{R}_{L}$, then for any $\bm     d\in
\mathcal{D}$ and for each $j=1,2,\ldots,m$, the long-run average profit
$\eta^{\bm     d}$ is strictly monotone decreasing with respect to each
decision element $d_{n,j}$, where $d_{n,j}\in\left\{  0,1,\ldots,j\right\}  $.
\end{The}

\textbf{Proof: }This proof is similar to the proof of Theorem \ref{The:left}.
For each $j=1,2,\ldots,m$, we consider two energy-efficient policies with an
interrelated structure as follows.%
\begin{align*}
\bm  d  &  =\left(  0,0,\ldots,0;d_{n,1},\ldots,d_{n,j-1}\underline{,d_{n,j}%
,}d_{n,j+1},\ldots,d_{n,m}\right)  ,\\
\bm  d^{\prime}  &  =\left(  0,0,\ldots,0;d_{n,1},\ldots,d_{n,j-1}%
\underline{,d_{n,j}^{\prime},}d_{n,j+1},\ldots,d_{n,m}\right)  ,
\end{align*}
where $d_{n,j}^{\prime}>d_{n,j}$, and $d_{n,j},d_{n,j}^{\prime}\in\left\{
0,1,\ldots,j\right\}  $. It is clear that%
\[
\eta^{\bm     d^{\prime}}-\eta^{\bm     d}=\mu_{2}\pi^{\left(  \bm
d^{\prime}\right)  }\left(  n,j\right)  \left(  d_{n,j}^{\prime}%
-d_{n,j}\right)  \left[  {G}^{\left(  \bm     d\right)  }\left(  n,j\right)
+{c}\right]  .
\]
If $0\leq R\leq{R}_{L}$, then it is seen from Proposition \ref{Pro:servicep}
that for any $\bm     d\in\mathcal{D}$ and for each $j=1,2,\ldots,m$,
${G}^{\left(  \bm     d\right)  }\left(  n,j\right)  +{c\leq0}$. Thus, we
obtain that for the two policies $\bm     d,\bm     d^{\prime}\in\mathcal{D}$
with $d_{n,j}^{\prime}>d_{n,j}$ and $d_{n,j},d_{n,j}^{\prime}\in\left\{
0,1,\ldots,j\right\}  $,%
\[
\eta^{\bm     d^{\prime}}<\eta^{\bm     d}.
\]
This shows that%
\[
\eta^{\bm     d}|_{d_{n,j}=1}>\eta^{\bm     d}|_{d_{n,j}=2}>\cdots>\eta^{\bm
d}|_{d_{n,j}=m-1}>\eta^{\bm     d}|_{d_{n,j}=m}.
\]
This completes the proof. \textbf{{\rule{0.08in}{0.08in}}}

When $0\leq R\leq{R}_{L}$, we use Figure 3 to provide an intuitive summary for
the main results given in Theorems~\ref{The:right} and \ref{Cor:left}.

\begin{figure}[th]
\centering          \includegraphics[width=12cm]{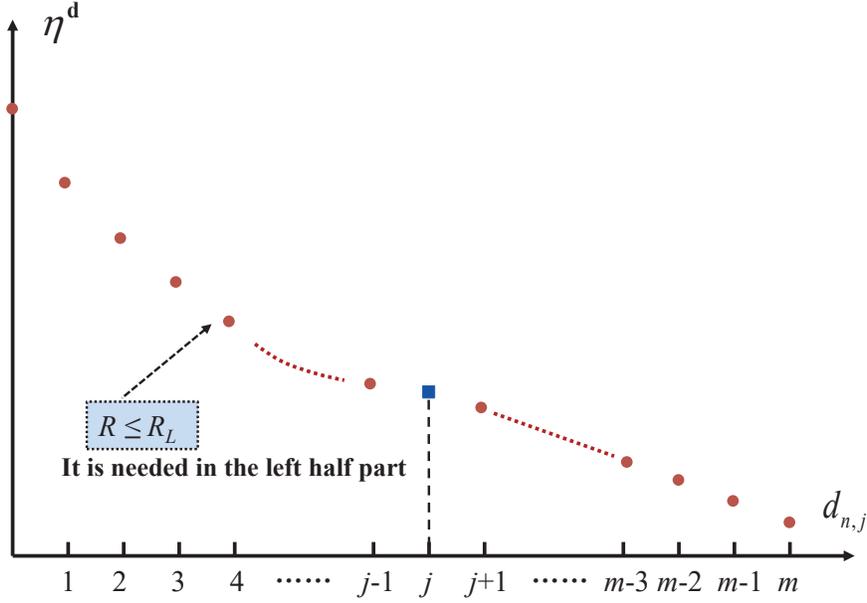}  \caption{The
decreasing structure of the long-run average profit.}%
\label{figure:fig-3}%
\end{figure}

The following theorem establishes the optimal energy-efficient policy $\bm
d^{\ast}$ in the data center, and also computes the maximal long-run average profit.

\begin{The}
The optimal energy-efficient policy $\bm     d^{\ast}$ and the maximal
long-run average profit $\eta^{\bm     d^{\ast}}$ can be determined in the
following two different cases:

(a) If $R\geq{R}_{H}$, then%
\[
\bm     d^{\ast}=\left(  0,0,\ldots,0;1,2,\ldots,m\right)
\]
and%
\begin{align*}
\eta^{\bm  d^{\ast}}  &  =\sum_{i=0}^{n}\frac{\dfrac{\lambda^{i}}{i!\mu
_{1}^{i}}}{\overset{n}{\underset{i=0}{\sum}}\dfrac{\lambda^{i}}{i!\mu_{1}^{i}%
}+\dfrac{\lambda^{n}}{n!\mu_{1}^{n}}\overset{m}{\underset{j=1}{\sum}}%
\dfrac{\lambda^{j}}{\underset{i=1}{\overset{j}{\Pi}}\left(  n\mu_{1}+i\mu
_{2}\right)  }}\left[  Ri\mu_{1}-\left(  nP_{1,W}+mP_{2,S}\right)
C_{1}-iC_{2}^{\left(  1\right)  }\right] \\
&  \text{ \ \ }+\sum_{j=1}^{m}\frac{\dfrac{\lambda^{n}}{n!\mu_{1}^{n}}%
\dfrac{\lambda^{j}}{\underset{i=1}{\overset{j}{\Pi}}\left(  n\mu_{1}+i\mu
_{2}\right)  }}{\overset{n}{\underset{i=0}{\sum}}\dfrac{\lambda^{i}}{i!\mu
_{1}^{i}}+\dfrac{\lambda^{n}}{n!\mu_{1}^{n}}\overset{m}{\underset{j=1}{\sum}%
}\dfrac{\lambda^{j}}{\underset{i=1}{\overset{j}{\Pi}}\left(  n\mu_{1}+i\mu
_{2}\right)  }}\Big\{  R\left(  n\mu_{1}+j\mu_{2}\right)  -\left[
nC_{2}^{\left(  1\right)  }+jC_{2}^{\left(  2\right)  }\right] \\
\text{ \ \ \ }  &  \text{\ \ \ }-\left[  nP_{1,W}+jP_{2,W}+\left(  m-j\right)
P_{2,S}\right]  C_{1}-n\mu_{1}1_{\{j>0\}}C_{3}-\lambda1_{\left\{
i=n,j=m\right\}  }C_{4}\Big\}  .
\end{align*}

(b) If $0\leq R\leq{R}_{L}$, then%
\[
\bm     d^{\ast}=\left(  0,0,\ldots,0;0,0,\ldots,0\right)
\]
and%
\begin{align*}
\eta^{\bm  d^{\ast}}  &  =\sum_{i=0}^{n}\frac{\dfrac{\lambda^{i}}{i!\mu
_{1}^{i}}}{\overset{n}{\underset{i=0}{\sum}}\dfrac{\lambda^{i}}{i!\mu_{1}^{i}%
}+\dfrac{\lambda^{n}}{n!\mu_{1}^{n}}\overset{m}{\underset{j=1}{\sum}}%
\dfrac{\lambda^{j}}{\left(  n\mu_{1}\right)  ^{j}}}\left[  Ri\mu_{1}-\left(
nP_{1,W}+mP_{2,S}\right)  C_{1}-iC_{2}^{\left(  1\right)  }\right] \\
&  \text{ \ \ }+\sum_{j=1}^{m}\frac{\dfrac{\lambda^{n}}{n!\mu_{1}^{n}}%
\dfrac{\lambda^{j}}{\left(  n\mu_{1}\right)  ^{j}}}{\overset{n}{\underset
{i=0}{\sum}}\dfrac{\lambda^{i}}{i!\mu_{1}^{i}}+\dfrac{\lambda^{n}}{n!\mu
_{1}^{n}}\overset{m}{\underset{j=1}{\sum}}\dfrac{\lambda^{j}}{\left(  n\mu
_{1}\right)  ^{j}}}\Big\{  Rn\mu_{1}-\left[  nP_{1,W}+mP_{2,S}\right]  C_{1}\\
\text{ \ \ }  &  \text{\ \ \ }-\left[  nC_{2}^{\left(  1\right)  }%
+jC_{2}^{\left(  2\right)  }\right]  -n\mu_{1}1_{\{j>0\}}C_{3}-\lambda
1_{\left\{  i=n,j=m\right\}  }C_{4}\Big\}  .
\end{align*}
\end{The}

\textbf{Proof: }(a) For the optimal energy-efficient policy $\bm     d^{\ast
}=\left(  0,0,\ldots,0;1,2,\ldots,m\right)  $, it is clear that $d_{n,j}%
^{\ast}=j$ and $d_{n,j}^{\ast}\wedge j=j$. Thus, it follows from
(\ref{eq-1-01}), (\ref{eq-1-02}) and (\ref{eq-1-4}) that
\[
\xi_{i,0}=\dfrac{\lambda^{i}}{i!\mu_{1}^{i}},\text{ \ }i=0,1,\ldots,n,
\]
and%
\[
\xi_{n,j}^{\left(  \bm     d^{\ast}\right)  }=\dfrac{\lambda^{n}}{n!\mu
_{1}^{n}}\dfrac{\lambda^{j}}{\underset{i=1}{\overset{j}{\Pi}}\left(  n\mu
_{1}+i\mu_{2}\right)  },\text{ \ }j=1,2,\ldots,m,
\]%
\[
b^{\left(  \bm     d^{\ast}\right)  }=\overset{n}{\underset{i=0}{\sum}}%
\dfrac{\lambda^{i}}{i!\mu_{1}^{i}}+\dfrac{\lambda^{n}}{n!\mu_{1}^{n}}%
\overset{m}{\underset{j=1}{\sum}}\dfrac{\lambda^{j}}{\underset{i=1}%
{\overset{j}{\Pi}}\left(  n\mu_{1}+i\mu_{2}\right)  }.
\]
It follows from (\ref{eq-1-2}) and (\ref{eq-1-3}) that for $i=0,1,\ldots,n,$%
\[
\pi^{\left(  \bm     d^{\ast}\right)  }\left(  i,0\right)  =\frac{\dfrac
{\lambda^{i}}{i!\mu_{1}^{i}}}{\overset{n}{\underset{i=0}{\sum}}\dfrac
{\lambda^{i}}{i!\mu_{1}^{i}}+\dfrac{\lambda^{n}}{n!\mu_{1}^{n}}\overset
{m}{\underset{j=1}{\sum}}\dfrac{\lambda^{j}}{\underset{i=1}{\overset{j}{\Pi}%
}\left(  n\mu_{1}+i\mu_{2}\right)  }},
\]
and for $j=1,2,\ldots,m,$%
\[
\pi^{\left(  \bm     d^{\ast}\right)  }\left(  n,j\right)  =\frac{\dfrac
{\lambda^{n}}{n!\mu_{1}^{n}}\dfrac{\lambda^{j}}{\underset{i=1}{\overset{j}%
{\Pi}}\left(  n\mu_{1}+i\mu_{2}\right)  }}{\overset{n}{\underset{i=0}{\sum}%
}\dfrac{\lambda^{i}}{i!\mu_{1}^{i}}+\dfrac{\lambda^{n}}{n!\mu_{1}^{n}}%
\overset{m}{\underset{j=1}{\sum}}\dfrac{\lambda^{j}}{\underset{i=1}%
{\overset{j}{\Pi}}\left(  n\mu_{1}+i\mu_{2}\right)  }}.\text{ }%
\]
At the same time, from (\ref{eq-50}) to (\ref{eq-54}) we obtain that for
$i=0,1,\ldots,n$,
\[
f\left(  i,0\right)  =Ri\mu_{1}-\left(  nP_{1,W}+mP_{2,S}\right)  C_{1}%
-iC_{2}^{\left(  1\right)  },
\]
and $j=1,2,\ldots,m$,
\begin{align*}
f^{\left(  \bm  d^{\ast}\right)  }\left(  n,j\right)   &  =R\left(  n\mu
_{1}+j\mu_{2}\right)  -\left[  nP_{1,W}+jP_{2,W}+\left(  m-j\right)
P_{2,S}\right]  C_{1}\\
&  \text{ \ \ }-\left[  nC_{2}^{\left(  1\right)  }+jC_{2}^{\left(  2\right)
}\right]  -n\mu_{1}1_{\{j>0\}}C_{3}-\lambda1_{\left\{  i=n,j=m\right\}  }%
C_{4}.
\end{align*}
Thus, we obtain%
\[
\eta^{\bm     d^{\ast}}=\sum_{i=0}^{n}\pi^{\left(  \bm     d^{\ast}\right)
}\left(  i,0\right)  f\left(  i,0\right)  +\sum_{j=1}^{m}\pi^{\left(  \bm
d^{\ast}\right)  }\left(  n,j\right)  f^{\left(  \bm     d^{\ast}\right)
}\left(  n,j\right)  .
\]
A simple computation directly derives our desired result.

(b) For the optimal policy $\bm         d^{\ast}=\left(  0,0,\ldots
,0;0,0,\ldots,0\right)  $, it is clear that $d_{n,j}^{\ast}=0$ so that
$d_{n,j}^{\ast}\wedge j=0$. A similar analysis to that in (a) can lead to our
desired result. This completes the proof. \textbf{{\rule{0.08in}{0.08in}}}

\begin{Rem}
The results of Theorem 6 are intuitive due to the fact that when the service
price is suitably high, the number of working servers is equal to the number
of waiting jobs in Group 2; while when the service price is lower, each server
opened at the work state will pay a high energy consumption cost but receive a
low revenue, thus the profit cannot increase and all the servers in Group 2
would like to be at the sleep state.
\end{Rem}

When the price $R_{L}<R<R_{H}$, we can further derive the following theorem
about the monotonicity of $\eta^{\bm                      d}$ with respect to
the decision element $d_{n,j}$.

\begin{The}
\label{The:mono} If $R_{L}<R<R_{H}$, then the long-run average profit
$\eta^{\bm     d}$ is monotone (either increasing or decreasing) with respect
to the decision element $d_{n,j}$, where $j=1,2,\dots,m$ and $d_{n,j}%
\in\left\{  0,1,\dots,j\right\}  $.
\end{The}

\textbf{Proof: }Similar to the first part of the proof for
Theorem~\ref{The:left}, we consider any two energy-efficient policies with an
interrelated structure as follows.%
\begin{align*}
\bm  d  &  =\left(  0,0,\ldots,0;d_{n,1},\ldots,d_{n,j-1}\underline{,d_{n,j}%
,}d_{n,j+1},\ldots,d_{n,m}\right)  ,\\
\bm  d^{\prime}  &  =\left(  0,0,\ldots,0;d_{n,1},\ldots,d_{n,j-1}%
\underline{,d_{n,j}^{\prime},}d_{n,j+1},\ldots,d_{n,m}\right)  ,
\end{align*}
where $d_{n,j},d_{n,j}^{\prime}\in\left\{  0,1,\ldots,j\right\}  $. Applying
Lemma~\ref{Lem:diff}, we obtain
\begin{equation}
\eta^{\bm  d^{\prime}}-\eta^{\bm  d}=\mu_{2}\pi^{\left(  \bm  d^{\prime
}\right)  }\left(  n,j\right)  \left(  d_{n,j}^{\prime}-d_{n,j}\right)
\left[  {G}^{\left(  \bm  d\right)  }\left(  n,j\right)  +{c}\right]  .
\label{eq-a1}%
\end{equation}
On the other hand, we can similarly obtain the following difference equation
\begin{equation}
\eta^{\bm  d}-\eta^{\bm  d^{\prime}}=\mu_{2}\pi^{\left(  \bm  d\right)
}\left(  n,j\right)  \left(  d_{n,j}-d_{n,j}^{\prime}\right)  \left[
{G}^{\left(  \bm  d^{\prime}\right)  }\left(  n,j\right)  +{c}\right]  .
\label{eq-a2}%
\end{equation}
By summing (\ref{eq-a1}) and (\ref{eq-a2}), we have
\[
\mu_{2}\pi^{\left(  \bm   d^{\prime}\right)  }\left(  n,j\right)  \left(
d_{n,j}^{\prime}-d_{n,j}\right)  \left[  {G}^{\left(  \bm   d\right)  }\left(
n,j\right)  +{c}\right]  +\mu_{2}\pi^{\left(  \bm   d\right)  }\left(
n,j\right)  \left(  d_{n,j}-d_{n,j}^{\prime}\right)  \left[  {G}^{\left(  \bm
d^{\prime}\right)  }\left(  n,j\right)  +{c}\right]  =0.
\]
We can directly derive
\[
\pi^{\left(  \bm   d^{\prime}\right)  }\left(  n,j\right)  \left[
{G}^{\left(  \bm   d\right)  }\left(  n,j\right)  +{c}\right]  =\pi^{\left(
\bm   d\right)  }\left(  n,j\right)  \left[  {G}^{\left(  \bm   d^{\prime
}\right)  }\left(  n,j\right)  +{c}\right]  .
\]
Therefore, we have the \emph{sign conservation equation}
\begin{equation}
\frac{{G}^{\left(  \bm  d\right)  }\left(  n,j\right)  +{c}}{{G}^{\left(  \bm
d^{\prime}\right)  }\left(  n,j\right)  +{c}}=\frac{\pi^{\left(  \bm
d\right)  }\left(  n,j\right)  }{\pi^{\left(  \bm  d^{\prime}\right)  }\left(
n,j\right)  }>0. \label{eq-sign-cons}%
\end{equation}
The above equation means that the sign of ${G}^{\left(  \bm   d\right)
}\left(  n,j\right)  +{c}$ and ${G}^{\left(  \bm   d^{\prime}\right)  }\left(
n,j\right)  +{c}$ are always identical when a particular decision element
$d_{n,j}$ is changed to any $d_{n,j}^{\prime}$. With the sign conservation
equation (\ref{eq-sign-cons}) and the performance difference equation
(\ref{eq-a2}), we can directly derive that the long-run average profit
$\eta^{\bm   d}$ is monotone with respect to $d_{n,j}$. This completes the
proof. \textbf{{\rule{0.08in}{0.08in}}}

Based on Theorems~\ref{The:left}, \ref{Cor:left}, and \ref{The:mono}, we can
directly derive that the optimal decision element $d_{n,j}^{\ast}$ has the
bang-bang control form, no matter what the value of $R$ will be.

\begin{Cor}
\label{Cor:bang-bang} The optimal decision element $d^{*}_{n,j}$ is either 0
or $j$, i.e., the bang-bang control is optimal.
\end{Cor}

With Corollary~\ref{Cor:bang-bang}, we should either keep all servers sleep or
turn on the servers such that the number of working servers equals the number
of waiting jobs in Group 2. We can see that the search space of $d_{n,j}$ can
be reduced from $\{0,1,\dots,j\}$ to a 2-element set $\{0,j\}$, which is a
significant reduction of optimization complexity. The form of the bang-bang
control is also very simple and it is easy to adopt in practice, while the
optimality of the bang-bang control guarantees the performance confidence of
such simple forms of control.

\section{Threshold Energy-Efficient Policy}

We have proved the optimality of the bang-bang control, no matter what the
value of $R$ will be. In practice, threshold-type policy is another category
of policies which also have a very simple form and are widely adopted in many
practical systems. In this section, we focus our study on the threshold-type
policy, although its optimality is not yet proved rigorously in our problem.
We use the Poisson equation to study a class of threshold energy-efficient
policies, and obtain the necessary condition of the optimal threshold
energy-efficient policy.

Now, we introduce an interesting subset of the policy set $\mathcal{D}$ as
follows. To this end, for $\theta=1,2,\ldots,m+1$, we write $\bm
d_{\theta}$ as an energy-efficient policy $\bm     d$ with $d_{n,j}=0$ if
$1\leq j\leq\theta-1\ $and $d_{n,j}=j$ if $\theta\leq j\leq m$, i.e.,
\[
\bm     d_{\theta}\overset{\text{def}}{=}\left(  0,0,\ldots,0;\underset
{\theta-1\text{ zeros}}{\underbrace{0,0,\ldots,0},}\theta,\theta
+1,\ldots,m\right)  .
\]

Let%
\[
\mathcal{D}^{\Delta}\overset{\text{def}}{=}\left\{  \bm    d_{\theta}%
:\theta=1,2,\ldots,m+1\right\}  .
\]
Then,%
\[
\mathcal{D}^{\Delta}=\left\{  \left(  0,0,\ldots,0;\underset{\theta-1\text{
zeros}}{\underbrace{0,0,\ldots,0},}\theta,\theta+1,\ldots,m\right)
:\theta=1,2,\ldots,m+1\right\}  .
\]
It is easy to see that $\mathcal{D}^{\Delta}\subset\mathcal{D}$.

For a policy $\bm   d_{\theta}$, it is clear that if $1\leq j\leq\theta-1$,
then $d_{n,j}=0$ and $d_{n,j}\wedge j=0$; and if $\theta\leq j\leq m$, then
$d_{n,j}=j$ and $d_{n,j}\wedge j=j$. Thus, it follows from (\ref{eq-1-01}),
(\ref{eq-1-02}) and (\ref{eq-1-4}) that
\[
\xi_{i,0}=\dfrac{\lambda^{i}}{i!\mu_{1}^{i}},\text{ \ }i=0,1,\ldots,n,
\]
and%
\[
\xi_{n,j}^{\left(  \bm   d_{\theta}\right)  }=\dfrac{\lambda^{n}}{n!\mu
_{1}^{n}}\dfrac{\lambda^{j}}{\left(  n\mu_{1}\right)  ^{j}},\text{
\ }j=1,2,\ldots,\theta-1,
\]%
\[
\xi_{n,j}^{\left(  \bm   d_{\theta}\right)  }=\dfrac{\lambda^{n}}{n!\mu
_{1}^{n}}\dfrac{\lambda^{\theta-1}}{\left(  n\mu_{1}\right)  ^{\theta-1}%
}\dfrac{\lambda^{j-\theta+1}}{\underset{i=\theta}{\overset{j}{\Pi}}\left(
n\mu_{1}+i\mu_{2}\right)  },\text{ \ }j=\theta,\theta+1,\ldots,m,
\]%
\[
b^{\left(  \bm   d_{\theta}\right)  }=\overset{n}{\underset{i=0}{\sum}}%
\dfrac{\lambda^{i}}{i!\mu_{1}^{i}}+\dfrac{\lambda^{n}}{n!\mu_{1}^{n}}%
\sum_{j=1}^{\theta-1}\dfrac{\lambda^{j}}{\left(  n\mu_{1}\right)  ^{j}}%
+\dfrac{\lambda^{n}}{n!\mu_{1}^{n}}\dfrac{\lambda^{\theta-1}}{\left(  n\mu
_{1}\right)  ^{\theta-1}}\sum_{j=\theta}^{m}\dfrac{\lambda^{j-\theta+1}%
}{\underset{i=\theta}{\overset{j}{\Pi}}\left(  n\mu_{1}+i\mu_{2}\right)  }.
\]
It follows from (\ref{eq-1-2}) and (\ref{eq-1-3}) that for $i=0,1,\ldots,n,$%
\[
\pi^{\left(  \bm   d_{\theta}\right)  }\left(  i,0\right)  =\frac{\dfrac
{\lambda^{i}}{i!\mu_{1}^{i}}}{\overset{n}{\underset{i=0}{\sum}}\dfrac
{\lambda^{i}}{i!\mu_{1}^{i}}+\dfrac{\lambda^{n}}{n!\mu_{1}^{n}}\sum
\limits_{j=1}^{\theta-1}\dfrac{\lambda^{j}}{\left(  n\mu_{1}\right)  ^{j}%
}+\dfrac{\lambda^{n}}{n!\mu_{1}^{n}}\dfrac{\lambda^{\theta-1}}{\left(
n\mu_{1}\right)  ^{\theta-1}}\sum\limits_{j=\theta}^{m}\dfrac{\lambda
^{j-\theta+1}}{\underset{i=\theta}{\overset{j}{\Pi}}\left(  n\mu_{1}+i\mu
_{2}\right)  }};
\]
for $j=1,2,\ldots,\theta-1,$%
\[
\pi^{\left(  \bm   d_{\theta}\right)  }\left(  n,j\right)  =\frac{\dfrac
{\lambda^{n}}{n!\mu_{1}^{n}}\dfrac{\lambda^{j}}{\left(  n\mu_{1}\right)  ^{j}%
}}{\overset{n}{\underset{i=0}{\sum}}\dfrac{\lambda^{i}}{i!\mu_{1}^{i}}%
+\dfrac{\lambda^{n}}{n!\mu_{1}^{n}}\sum\limits_{j=1}^{\theta-1}\dfrac
{\lambda^{j}}{\left(  n\mu_{1}\right)  ^{j}}+\dfrac{\lambda^{n}}{n!\mu_{1}%
^{n}}\dfrac{\lambda^{\theta-1}}{\left(  n\mu_{1}\right)  ^{\theta-1}}%
\sum\limits_{j=\theta}^{m}\dfrac{\lambda^{j-\theta+1}}{\underset{i=\theta
}{\overset{j}{\Pi}}\left(  n\mu_{1}+i\mu_{2}\right)  }};
\]
and for $j=\theta,\theta+1,\ldots,m,$%
\[
\pi^{\left(  \bm   d_{\theta}\right)  }\left(  n,j\right)  =\frac{\dfrac
{\lambda^{n}}{n!\mu_{1}^{n}}\dfrac{\lambda^{\theta-1}}{\left(  n\mu
_{1}\right)  ^{\theta-1}}\dfrac{\lambda^{j-\theta+1}}{\underset{i=\theta
}{\overset{j}{\Pi}}\left(  n\mu_{1}+i\mu_{2}\right)  }}{\overset{n}%
{\underset{i=0}{\sum}}\dfrac{\lambda^{i}}{i!\mu_{1}^{i}}+\dfrac{\lambda^{n}%
}{n!\mu_{1}^{n}}\sum\limits_{j=1}^{\theta-1}\dfrac{\lambda^{j}}{\left(
n\mu_{1}\right)  ^{j}}+\dfrac{\lambda^{n}}{n!\mu_{1}^{n}}\dfrac{\lambda
^{\theta-1}}{\left(  n\mu_{1}\right)  ^{\theta-1}}\sum\limits_{j=\theta}%
^{m}\dfrac{\lambda^{j-\theta+1}}{\underset{i=\theta}{\overset{j}{\Pi}}\left(
n\mu_{1}+i\mu_{2}\right)  }}.\text{ }%
\]
It follows from (\ref{eq-50}) to (\ref{eq-54}) that for $i=0,1,\ldots,n$,
\[
f\left(  i,0\right)  =Ri\mu_{1}-\left(  nP_{1,W}+mP_{2,S}\right)  C_{1}%
-iC_{2}^{\left(  1\right)  };
\]
for $j=1,2,\ldots,\theta-1$,
\[
f^{\left(  \bm   d_{\theta}\right)  }\left(  n,j\right)  =Rn\mu_{1}-\left(
nP_{1,W}+mP_{2,S}\right)  C_{1}-\left[  nC_{2}^{\left(  1\right)  }%
+jC_{2}^{\left(  2\right)  }\right]  -n\mu_{1}C_{3};
\]
and for $j=\theta,\theta+1,\ldots,m,$
\begin{align*}
f^{\left(  \bm  d_{\theta}\right)  }\left(  n,j\right)   &  =R\left(  n\mu
_{1}+j\mu_{2}\right)  -\left[  nP_{1,W}+jP_{2,W}+\left(  m-j\right)
P_{2,S}\right]  C_{1}\\
&  -\left[  nC_{2}^{\left(  1\right)  }+jC_{2}^{\left(  2\right)  }\right]
-n\mu_{1}C_{3}-\lambda1_{\left\{  i=n,j=m\right\}  }C_{4}.
\end{align*}
Note that%
\[
\eta^{\bm   d_{\theta}}=\sum_{i=0}^{n}\pi^{\left(  \bm   d_{\theta}\right)
}\left(  i,0\right)  f\left(  i,0\right)  +\sum_{j=1}^{m}\pi^{\left(  \bm
d_{\theta}\right)  }\left(  n,j\right)  f^{\left(  \bm   d_{\theta}\right)
}\left(  n,j\right)  .
\]
We obtain the explicit expression of the long-run average profit under policy
$\bm   d_{\theta}$ as follows.%
\begin{align*}
\eta^{\bm  d_{\theta}}  &  =\sum_{i=0}^{n}\pi^{\left(  \bm  d_{\theta}\right)
}\left(  i,0\right)  f\left(  i,0\right)  \left[  Ri\mu_{1}-\left(
nP_{1,W}+mP_{2,S}\right)  C_{1}-iC_{2}^{\left(  1\right)  }\right] \\
&  \text{ \ \ \ }+\sum_{j=1}^{\theta-1}\pi^{\left(  \bm  d_{\theta}\right)
}\left(  n,j\right)  \left\{  Rn\mu_{1}-\left(  nP_{1,W}+mP_{2,S}\right)
C_{1}-\left[  nC_{2}^{\left(  1\right)  }+jC_{2}^{\left(  2\right)  }\right]
-n\mu_{1}C_{3}\right\} \\
&  \text{ \ \ \ }+\sum_{j=\theta}^{m}\pi^{\left(  \bm  d_{\theta}\right)
}\left(  n,j\right)  \Big\{  R\left(  n\mu_{1}+j\mu_{2}\right)  -\left[
nP_{1,W}+jP_{2,W}+\left(  m-j\right)  P_{2,S}\right]  C_{1}\\
&  \text{ \ \ \ }-\left[  nC_{2}^{\left(  1\right)  }+jC_{2}^{\left(
2\right)  }\right]  -n\mu_{1}C_{3}-\lambda1_{\left\{  i=n,j=m\right\}  }%
C_{4}\Big\}  .\text{\ \ \ \ }%
\end{align*}

Let%
\[
\theta^{\ast}=\underset{\theta\in\left\{  1,2,\ldots,m+1\right\}  }{\arg\max
}\left\{  \eta^{\bm     d_{^{\theta}}}\right\}  .
\]
Then, we call $\bm     d_{\theta^{\ast}}$ the optimal threshold
energy-efficient policy in the policy set $\mathcal{D}^{\Delta}$. Since
$\mathcal{D}^{\Delta}\subset\mathcal{D}$, the partially ordered set
$\mathcal{D}$ shows that $\mathcal{D}^{\Delta}$ is also a partially ordered
set. Based on this, it is easy to see from the two partially ordered sets
$\mathcal{D}$ and $\mathcal{D}^{\Delta}$ that%
\[
\eta^{\bm     d_{\theta^{\ast}}}\leq\eta^{\bm     d^{\ast}}.
\]
If $\eta^{\bm     d_{\theta^{\ast}}}=\eta^{\bm     d^{\ast}}$, then we call
$\bm   d_{\theta^{\ast}}$ the optimal threshold energy-efficient policy in the
original policy set $\mathcal{D}$; If $\eta^{\bm     d_{\theta^{\ast}}}%
<\eta^{\bm     d^{\ast}}$, then we call $\bm     d_{\theta^{\ast}}$ the
suboptimal threshold energy-efficient policy in the original\ policy set
$\mathcal{D}$.

We take a minimal positive integer $\theta^{\ast}\in\left\{  1,2,\ldots
,m+1\right\}  $ such that%
\[
\bm             d_{\theta^{\ast}}=\left(  0,0,\ldots,0;\underset{\theta^{\ast
}-1\text{ zeros}}{\underbrace{0,0,\ldots,0},}\theta^{\ast},\theta^{\ast
}+1,\ldots,m\right)  .
\]

For the optimal threshold energy-efficient policy $\bm     d_{\theta^{\ast}}$,
the following theorem determines the positive or negative property of the
function $G^{\left(  \bm     d_{\theta}\right)  }\left(  n,\theta\right)  +c$
for $\theta=\theta^{\ast}-1,\theta^{\ast},\theta^{\ast}+1$, although the
explicit expression of the perturbation realization factor $G^{\left(  \bm
d_{\theta}\right)  }\left(  n,\theta\right)  $ is not given yet. This may be
useful for us to understand the role played by Proposition 2 in analyzing the
monotonicity and optimality of the energy-efficient policies. Furthermore, we
also derive the necessary condition of the optimal threshold energy-efficient policy.

\begin{The}
For the threshold energy-efficient policies of the data center, the optimal
threshold policy $\bm            d^{*}_{\theta}$ satisfies the following
condition
\[
G^{\left(  \bm                      d_{\theta^{\ast}-1}\right)  }\left(
n,\theta^{\ast}-1\right)  +{c\leq0,}\text{ \ \ }G^{\left(  \bm  d_{\theta
^{\ast}}\right)  }\left(  n,\theta^{\ast}\right)  +{c\geq0,}\text{ \ \ }%
{G}^{\left(  \bm                  d_{\theta^{\ast}+1}\right)  }\left(
n,\theta^{\ast}+1\right)  +{c\geq0}.
\]
\end{The}

\textbf{Proof:} We consider three threshold energy-efficient policies with an
interrelated structure as follows.%
\begin{align*}
\bm  d_{\theta^{\ast}}  &  =\left(  0,0,\ldots,0;\underset{\mathbf{\theta
}^{\ast}-1\text{ zeros}}{\underbrace{0,0,\ldots,0},}\theta^{\ast},\theta
^{\ast}+1,\ldots,m\right)  ,\\
\bm  d_{\theta^{\ast}+1}  &  =\left(  0,0,\ldots,0;\underset{\mathbf{\theta
}^{\ast}-1\text{ zeros}}{\underbrace{0,0,\ldots,0},}0,\theta^{\ast}%
+1,\ldots,m\right)  ,\\
\bm  d_{\theta^{\ast}-1}  &  =\left(  0,0,\ldots,0;\underset{\mathbf{\theta
}^{\ast}-2\text{ zeros}}{\underbrace{0,0,\ldots,0},}\theta^{\ast}%
-1,\theta^{\ast},\ldots,m\right)  .
\end{align*}
It follows from Lemma 1 that for two energy-efficient policies with an
interrelated structure
\begin{align*}
\bm  d  &  =\left(  0,0,\ldots,0;d_{n,1},\ldots,d_{n,j-1}\underline{,d_{n,j}%
,}d_{n,j+1},\ldots,d_{n,m}\right)  ,\\
\bm  d^{\prime}  &  =\left(  0,0,\ldots,0;d_{n,1},\ldots,d_{n,j-1}%
\underline{,d_{n,j}^{\prime},}d_{n,j+1},\ldots,d_{n,m}\right)  ,
\end{align*}
it is clear that%
\[
\eta^{\bm     d^{\prime}}-\eta^{\bm     d}=\mu_{2}\pi^{\left(  \bm
d^{\prime}\right)  }\left(  d_{n,j}^{\prime}-d_{n,j}\right)  \left[
G^{\left(  \bm  d\right)  }\left(  n,j\right)  +c\right]  .
\]
Thus, we obtain%
\[
\eta^{\bm     d_{\theta^{\ast}+1}}-\eta^{\bm     d_{\theta^{\ast}}}%
=-\theta^{\ast}\mu_{2}\pi^{\left(  \bm     d_{\theta^{\ast}+1}\right)
}\left[  G^{\left(  \bm    d_{\theta^{\ast}}\right)  }\left(  n,\theta^{\ast
}\right)  +c\right]  ,
\]
which, together with $\eta^{\bm     d_{\theta^{\ast}+1}}-\eta^{\bm
d_{\theta^{\ast}}}\leq0$, leads to%
\[
G^{\left(  \bm     d_{\theta^{\ast}}\right)  }\left(  n,\theta^{\ast}\right)
+{c\geq0.}%
\]
Similarly, we have%
\[
\eta^{\bm     d_{\theta^{\ast}}}-\eta^{\bm     d_{\theta^{\ast}+1}}%
=\theta^{\ast}\mu_{2}\pi^{\left(  \bm     d_{\theta^{\ast}}\right)  }\left[
{G}^{\left(  \bm    d_{\theta^{\ast}+1}\right)  }\left(  n,\theta^{\ast
}+1\right)  +c\right]  ,
\]
which indicates%
\[
{G}^{\left(  \bm     d_{\theta^{\ast}+1}\right)  }\left(  n,\theta^{\ast
}+1\right)  +{c\geq0.}%
\]
Also, we have%
\[
\eta^{\bm     d_{\theta^{\ast}}}-\eta^{\bm     d_{\theta^{\ast}-1}}=-\left(
\theta^{\ast}-1\right)  \mu_{2}\pi^{\left(  \bm     d_{\theta^{\ast}}\right)
}\left[  {G}^{\left(  \bm     d_{\theta^{\ast}-1}\right)  }\left(
n,\theta^{\ast}-1\right)  +{c}\right]  ,
\]
which indicates%
\[
{G}^{\left(  \bm     d_{\theta^{\ast}-1}\right)  }\left(  n,\theta^{\ast
}-1\right)  +{c\leq0.}%
\]
This completes the proof. \textbf{{\rule{0.08in}{0.08in}}}

\section{Conclusion}

In this paper, we propose a novel dynamic decision method by applying the
sensitivity-based optimization theory to study the optimal energy-efficient
policy of a data center with two groups of heterogeneous servers. We propose a
job transfer rule among the group-servers such that the sleep energy-efficient
mechanism of Group 2 becomes more effective. To find the optimal
energy-efficient policy of the data center, we set up a policy-based Poisson
equation and provide explicit expression for its unique solution by means of
the RG-factorization. Based on this, we derive the monotonicity and optimality
of the long-run average profit with respect to the energy-efficient policies
under some restrained service prices. We prove the optimality of the bang-bang
control, which significantly reduces the action search space. We also study
the threshold energy-efficient policy and derive the necessary condition of
the optimal threshold policy. Different from previous works in the literature
on applying the traditional MDP theory to the dynamic control of data centers,
the sensitivity-based optimization method used in this paper is easier and
more convenient in the study of energy-efficient data centers. This
sensitivity-based optimization method may open a new avenue to study the
optimal energy-efficient policy for more complicated data centers.

Along such a research line of applying the sensitivity-based optimization and
the RG-factorization to the energy-efficient data centers, the extension to
multiple groups of heterogeneous servers deserves further investigations. The
control policy will become more complicated when multiple groups of servers
are considered. Another interesting research topic is to consider different
cost structures, waiting capacities, service disciplines, or job migration
rules. Especially, when the job migration is not allowed in data centers, the
complexity of the dynamic control problem will dramatically increase and it
deserves further more investigations.

\section*{Acknowledgements}

Li Xia was supported by the National Key Research and Development Program of
China (2016YFB0901900, 2017YFC0704100), the National Natural Science
Foundation of China under grant No. 61573206 and No. U1301254, the National
111 International Collaboration Project (B06002), and the Suzhou-Tsinghua
Innovation Leading Action Project.

Quan-Lin Li was supported by the National Natural Science Foundation of China
under grant No. 71671158 and No. 71471160, and by the Natural Science
Foundation of Hebei province under grant No. G2017203277.

\end{document}